\documentclass[12pt, a4paper]{article} 

\usepackage[utf8]{inputenc} 

\usepackage[english]{babel}

\usepackage{graphicx}
\usepackage{color}
\usepackage{mathtools}
\usepackage{amssymb, amsmath}

\usepackage[percent]{overpic}

\graphicspath{ {./} }

\newtheorem{definition}{Definition}
\newtheorem{theorem}{Theorem}
\newtheorem{lemma}{Lemma}

\newtheorem{corollary}{Corollary}
\newtheorem{remark}{Remark}

\DeclareMathOperator{\sgn}{sgn}

\usepackage{subcaption}

\usepackage{float}

\title{The V-line transform with some generalizations and cone differentiation}

\author{G. Ambartsoumian$^\dagger$ {\small and} M. J. Latifi Jebelli$^*$\\
\small $^\dagger$Department of Mathematics, University of Texas at Arlington\\
\small $^*$Department of Mathematics, University of Arizona}

\date{}

\begin{document}

\maketitle


\begin{abstract}
The paper studies various properties of the V-line transform (VLT) in the plane and the conical Radon transform (CRT) in $\mathbb{R}^n$. The VLT maps a function to a family of its integrals along trajectories made of two rays emanating from a common point. The CRT considered in this paper maps a function to a set of its integrals over surfaces of polyhedral cones. These types of operators appear in mathematical models of single scattering optical tomography, Compton camera imaging and other applications. We derive new explicit inversion formulae for the VLT and the CRT, as well as proving some previously known results using more intuitive geometric ideas. Using our inversion formula for the VLT, we describe the range of that transformation when applied to a fairly broad class of functions and prove some support theorems. The efficiency of our method is demonstrated on several numerical examples. As an auxiliary result that plays a big role in this article, we derive a generalization of the Fundamental Theorem of Calculus, which we call the Cone Differentiation Theorem.
\end{abstract}


\section{Introduction}

The V-line transform is a generalized Radon transform, which puts into correspondence to a given function its integrals along ``broken lines'', i.e. piecewise linear trajectories that consist of two rays emanating from one vertex. The name of the operator is due to the resemblance of its integration trajectories to the letter ``V''. There is also a closely related broken ray transform (BRT), which integrates functions along broken rays (i.e. one of the branches of ``V'' has a finite length). If these transforms are considered on functions with compact support and the origin of the broken ray is outside of the convex hull of the support, then there is essentially no difference between the BRT and the VLT.

The VLT and its generalizations have attracted significant interest from the mathematical community in recent years. Part of this interest can be attributed to the development of various imaging techniques that use these transforms as a basis of their mathematical models. But in many other cases the study of such operators is of purely mathematical interest with intriguing connections between integral geometry, harmonic analysis, PDE's, microlocal analysis, differential geometry and other areas of mathematics.

To motivate the study of the VLT and its generalizations, we start with a brief description of two imaging modalities that use such transformations: single scattering optical tomography and Compton scattering tomography.

\begin{figure}[h]
\begin{center}
\includegraphics[width=50mm,keepaspectratio]{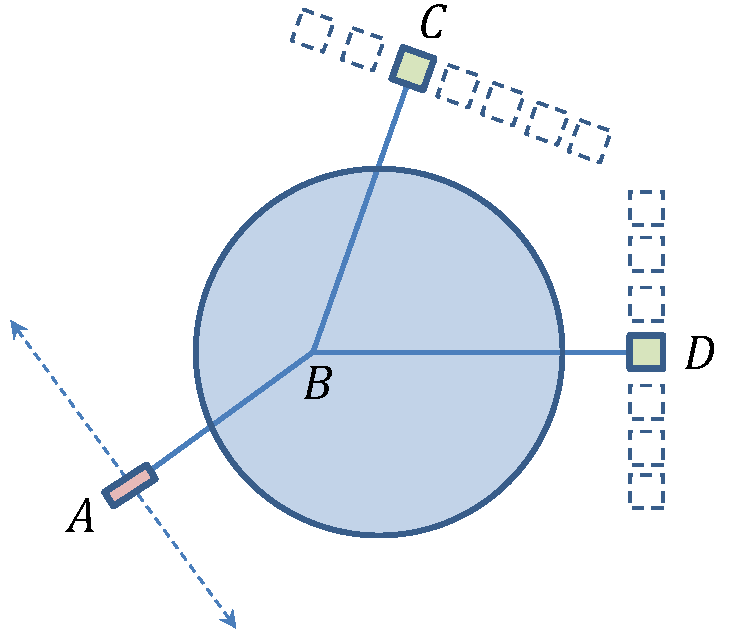}
\end{center}
\caption{Photons emitted in a fixed direction by the source $A$ are scattered inside the body. Arrays of collimated receivers catch the photons scattered in the directions of their collimation. The data is
collected for multiple positions of $A$ along a line perpendicular to the initial direction of photon beams. Knowing the locations of the source $A$ and the collimated receiver (e.g. $C$ or $D$), one can uniquely recover the scattering location $B$. Then, matching the VLT data measured by receivers $C$ and $D$ for the same source $A$ and scattering location $B$, one can generate the signed VLT data corresponding to the V-line $BCD$.}
\label{fig:SSOT}
\end{figure}

Optical tomography uses measurements of light transmitted through or scattered inside a biological object to recover the internal structure of that object. If the object is optically thin, then the majority of light photons fly through the object without scattering, and the measurements of outgoing light can be modeled by the (ordinary) Radon transform. If the object is thick, then the majority of photons go through multiple scattering events inside the object, and the corresponding process is usually modeled by the diffusion equation. In the case when the object is of certain moderate thickness, one can assume that the light photons scatter at most once inside the object (see Figure \ref{fig:SSOT}). Image reconstruction from such measurements is called single scattering optical tomography (SSOT). The measured data set in SSOT corresponds to a generalized Radon transform integrating the light attenuation coefficient along broken rays that coincide with the trajectories of scattered photons \cite{FMS-PhysRev-10, Florescu-Markel-Schotland,  FMS-PhysRev-09}. Hence, one of the crucial mathematical tasks in SSOT is the inversion of the BRT (or the VLT) in various geometric setups.

Notice that, by subtracting the measured data corresponding to the same source (e.g. $A$), the same scattering point (e.g. $B$), but different receivers (e.g. $C$ and $D$), one obtains an integral of the light attenuation coefficient along a V-line (in our case $CBD$) in which integration along each ray is done with a different algebraic sign (see Figure \ref{fig:SSOT}). Such ``signed'' VLT's have been studied before in \cite{Florescu-Markel-Schotland, Kats_Krylov-13}.

A slightly more complicated mathematical model is used in  tomographic reconstructions using single-scattered x-rays and accounting for energy dependent x-ray attenuation \cite{Krylov_Kats-15}. In this modality the x-ray photons are scattered by the charged particles of the matter. The scattering angle depends on the amount of energy lost by the photons. Hence, by using collimated receivers, one can register scattered x-ray photons that have lost the same amount of energy. Given the fact that the x-ray attenuation depends on the energy level of the x-ray, its integrals before and after the scattering are computed with different weights. As a result, here the image reconstruction problem requires inversion of a weighted VLT, where the unknown function is integrated with a different weight along each ray of the V-line.

It is easy to notice that the problem of inverting the VLT from full data is over-determined. The family of V-lines in the plane is 4-dimensional, while the image function depends on 2 variables. Similarly in 3D the family of V-lines is 7-dimensional, while the image function depends only on 3 variables. There are multiple options of limiting the set of V-lines to a subset of the appropriate dimension, e.g. limiting the locations of vertices, fixing the opening angles, fixing or limiting the axes of symmetry, etc. The choice of the appropriate setup for study is usually made based on the application at hand, as well as the mathematical considerations, e.g. the possibility and level of difficulty of inverting the transform.

Notice, that the VLT arising in SSOT has the vertices of integration trajectories {\em inside} the support of the image function. This feature distinguishes the mathematical problems arising in SSOT from those in Compton camera imaging discussed later in this section.

The first inversion formula for the VLT with vertices inside the support of the image function was presented in \cite{FMS-PhysRev-10, Florescu-Markel-Schotland,  FMS-PhysRev-09}. Here the authors considered V-lines in 2D slab geometry with a fixed opening angle, fixed axis of symmetry and arbitrary location of the vertex.  Simpler inversion formulae for the same setup were obtained later using other approaches in \cite{Gouia_Amb_V-line, Kats_Krylov-13, Sherson}. The VLT in other 2D geometries was studied in \cite{Ambartsoumian-VLine, Ambartsoumian_Moon_broken_ray_article,  Gaik_Souvik, Kats_Krylov-13, Krylov_Kats-15}. Inversion formulae for a generalization of the VLT to higher dimensions, namely the Conical Radon Transform mapping a function to a family of circular cones with vertices inside the image domain, were presented in \cite{Gouia_Conical_nd, Gouia_Amb_V-line, Palamodov}.

\begin{figure}[h]
\begin{center}
\includegraphics[width=60mm,keepaspectratio]{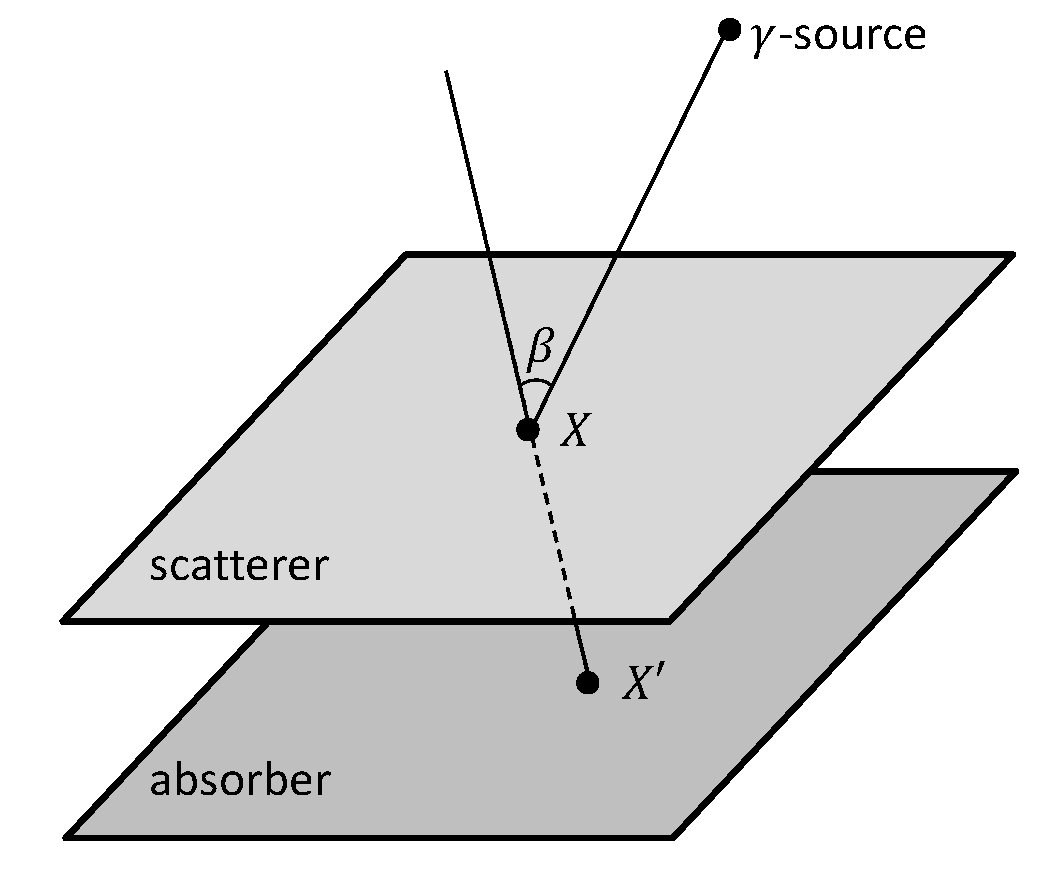}
\end{center}
\caption{A simple sketch of a Compton camera.}
\label{fig: Compton}
\end{figure}

Compton cameras are imaging devices that are used primarily for detection of sources of $\gamma$-radiation. These devices have a wide range of applications including astronomy, medicine and homeland security. A typical Compton camera consists of two parallel digital detectors: a scatterer and an absorber (see Figure \ref{fig: Compton}). When a photon in a $\gamma$-ray hits the scatterer at a point $X$, it changes its flight trajectory and hits the absorber at another point $X'$. The detectors of the Compton camera register these locations $X$ and $X'$, as well as the energy of the particle at each point. The well-known Compton scattering relation then allows recovery of the scattering angle $\beta$:
$$
\cos \beta = 1- \displaystyle\frac{mc^2\triangle E}{(E-\triangle E)E},
$$
where $E$ is the original energy of the photon, $\triangle E$ is its lost energy after scattering and $m$ is the mass of an electron.

Since the measured data set does not allow angularly resolving the location of the $\gamma$-source, one can only assume that the measurements correspond to a Radon-type transform over conical surfaces of the source distribution function. The mathematical task of image reconstruction here then corresponds to the inversion of such a CRT (e.g. see \cite{Allmaras, Basko_et_al, Cree_Bones, Maxim:FBP4cones, Maxim:Compton-full}). Note that in this case the vertices of cones of integration are limited to the surface of the scattering detector, and the image function ($\gamma$-source distribution) is supported on one side of that surface. It is easy to notice that even with such a restriction the problem of inverting the CRT is over-determined. The family of all cones with vertices located on a given surface is 5-dimensional, while the image function depends on 3 variables. There are multiple options of restricting the 5-dimensional set of cones to a 3-dimensional family, e.g. by using various combinations of fixing the direction of their axes of symmetry, fixing their opening angle, limiting the vertices to a curve, etc. One can also consider a 2D version of the same problem, where linear detectors are used instead of plane detectors. In that case the CRT becomes a VLT with vertices of broken rays on a line. Many researchers have obtained interesting results on the CRT and the VLT for such setups (e.g. see \cite{ Haltmeier-cones, Halt-Moon-Schief, Hristova-V-line, Jung_Moon-IP, Jung_Moon-SIAM,  Kuch_Terz_SIAM, Kuch_Terz_IPI,  Moon-conical-SIMA-16, Moon-conical-IP-17, Moon-Halt-Con_Cyl, MNTZ10, Nguyen_Truong_Grangeat, Schief-Halt_cone_sphere, Fatma_IP, Truong_Nguyen_V-line-1, Truong_Nguyen_V-line-2, TNZ07}). A nice survey of this field was recently published in \cite{TKK-review}.

Finally, we would like to mention a few other transforms that are related to the VLT and the BRT, but have significant differences. An operator, called the star transform, integrates a function along a star-like set that consists of multiple rays emanating from the same vertex (see \cite{ZSM-star_transform}). Another interesting area of research in integral geometry is dedicated to the recovery of functions defined inside a compact domain from their integrals along piecewise-linear trajectories that reflect multiple times from the boundary of that domain. As it often happens in mathematics, this transform is also called a broken-ray transform, although it is quite different from the BRT mentioned above. For more details and interesting results in this field we refer the reader to \cite{Hubenthal2014, Hubenthal2015, Ilmavirta-IP-2013, Ilmavirta-Salo-16}. In the setup of manifolds one can consider broken geodesics and various problems related to them (e.g. see \cite{Broken-geodesic}). Single scattering data are used in other imaging modalities besides SSOT and Compton cameras, e.g. in single-photon emission computed tomography (SPECT) \cite{SPECT-SS}.\\

In this article we consider the VLT and its high dimensional generalization integrating over polyhedral cones in cases with arbitrary vertex, but fixed axis of symmetry and opening angle.

The rest of the article is organized as follows. In Section \ref{Chap_cone} we derive a generalization of the Fundamental Theorem of Calculus from the perspective of partially ordered sets. That result, which we call the Cone Differentiation Theorem, is used throughout the rest of the paper as a major tool for studying the VLT and its generalizations. In Section \ref{Chap_BRT} we give formal definitions of the VLT with various weights and use the Cone Differentiation Theorem to derive both new and some previously known inversion formulae for the weighted VLT in the plane.  We then describe the range of the cone integration operator, and use it to characterize the range of VLT. Finally, we present here some support theorems using the corresponding theoretical discoveries from Section \ref{Chap_cone}. In Section \ref{Chap-nd} we generalize our inversion formula to the case of transforms integrating the image function along polyhedral cones in $\mathbb{R}^n$. In Section \ref{Chap_num} we present various numerical simulations of our inversion formulae. We list some additional remarks in Section \ref{Chap_remarks}, and summarize the results of the paper in Section \ref{Chap_sum}.  Some of the technical proofs have been separated as appendices in Sections \ref{Chap_ap1} and \ref{Chap_ap2}. 
\section{Cone Differentiation and Integration}\label{Chap_cone}

In this section we derive a generalization of the Fundamental Theorem of Calculus (FTC) to $\mathbb{R}^n$ from the perspective of partially ordered sets. We call this generalization the {\it Cone Differentiation Theorem} due to the concept of a positive cone in a partially ordered vector space.

We start our discussion by writing the FTC in terms of the natural order in $\mathbb{R}$ and then build the necessary background for our generalization to $\mathbb{R}^n$. If $\leq$ is the natural order on $\mathbb{R}$, for an integrable function $f$ and
\begin{equation}\label{FTC-1}
  F(x) = \int_{y\leq x} f(y)\,dy
\end{equation}
 we have $F^{'}=f$ almost everywhere.

\begin{remark}\label{abs_cont_1}
 Note that in this case $F$ is absolutely continuous.
\end{remark}
\subsection{Partial Order on $\mathbb{R}^n$}

Let us recall the concept of a partial order in a vector space.
A { \em partially ordered vector space} $V$ is a vector space over $\mathbb{R}$ together with a partial order $\leq$ such that:

\begin{enumerate}
\item If $x \leq y$, then $x+z \leq y+z$ for all $z \in V$.
\item If $x \geq 0$, then $cx \geq 0$ for all $c \in \mathbb{R}^{+}$.
\end{enumerate}

From the above definition we have $x \leq y \Leftrightarrow 0 \leq x-y$, and hence the order is completely determined by $V^+ = \{x \in V; x \geq 0\}$, which is called the {\em positive cone of $V$}.

For example, one can define a partial order in $\mathbb{R}^2$ as follows: $(x_1,x_2)\le(y_1,y_2)$ if and only if $x_1\le y_1$ and $x_2\le y_2$. In this case the positive cone will coincide with the first quadrant.

As another example, let us start with a choice of a positive cone in $\mathbb{R}^2$ and deduce the partial order from it. Consider $V^+=\{(x_1,x_2):\;x_1\ge 0, x_2\ge x_1\}$. Then the corresponding partial order in $\mathbb{R}^2$ will be given by the following: $(x_1,x_2)\le(y_1,y_2)$ if and only if $x_1\le y_1$ and $y_2-x_2\le y_1 - x_1$.

In general, for $P \subset V$ there is a partial order on $V$ such that $P=V^+$ if and only if
\[
P \cap (-P) = \{0\},
\]
\[
P+P \subset P,
\]
\[
c\geq 0 \Rightarrow cP \subset P.
\]

One can also identify a partial order structure with $V^-=-V^+$, called the {\em negative cone} of $V$.
It is easy to check that for the case of $V=\mathbb{R}^n$ the positive (and hence the negative) cone is actually a cone in the geometric sense.

These concepts play an important role in functional analysis and its applications.  For a detailed treatment of the subject we refer the reader to  \cite{Fremlin_book, Holmes}. In this paper we use the positive and negative cones to motivate and guide the generalization of certain classical results of analysis on the real line to higher dimensions.

Let us consider partial orders in $\mathbb{R}^n$ corresponding to positive (or negative) cones $C_{\mathcal{B}}$ generated by a set of fixed basis vectors $\mathcal{B}= \{v_1, ...,v_n\}$, i.e. $C_{\mathcal{B}} = \{ \sum_{i=1}^{n} c_iv_i; c_i \geq 0\}$.

 In the case of $\mathbb{R}^2$ we will use two linearly independent vectors $u,v$ as a generating set for the positive (or negative) cone. In this case the boundary of the cone is a V-line. This fact is an important building block of our construction of the inversion formula for the VLT.\\

\begin{figure}[h]
\begin{center}
\includegraphics[width=50mm,keepaspectratio]{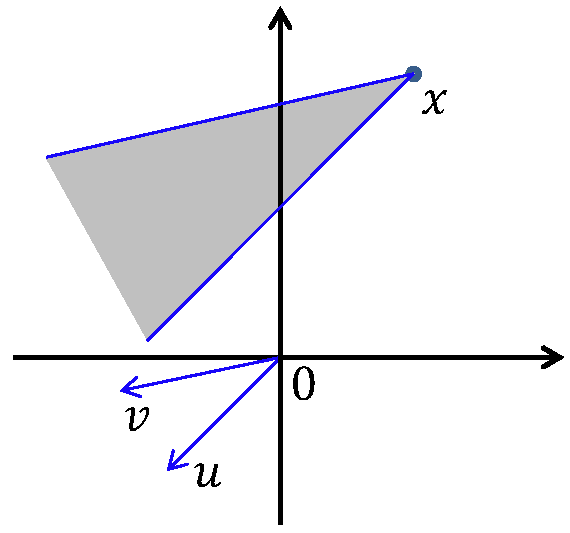}
\end{center}
\caption{The negative cone at $x$ with respect to a fixed partial order in $\mathbb{R}^2$ generated by vectors $u,v$.}
\label{fig: neg_cone}
\end{figure}

In analogy with formula (\ref{FTC-1}), for $f \in L^1(\mathbb{R}^n)$  we define $F$ on $\mathbb{R}^n$ as
\begin{equation}\label{FTC-n}
F(x) = \int_{y\leq x} f(y)\,d\mu,
\end{equation}
where $\mu$ is the standard Lebesgue measure on $\mathbb{R}^n$ and $y\leq x$ represents  the negative cone at $x$ with respect to a fixed partial order on $\mathbb{R}^n$. In other words, the integral is taken over the region $\{ y\in \mathbb{R}^n ; y\leq x \}$ (see Figure \ref{fig: neg_cone}).


\subsection{Two Classical Generalizations of FTC}

Now a natural question is: in what sense of differentiation can one generalize the FTC to $\mathbb{R}^n$?
While we construct a version of such a generalization using an order structure on $\mathbb{R}^n$ and its geometric properties, it should be mentioned that there are several classical results, which one can consider as generalizations of the FTC in a certain sense. In the derivation of our results we rely on two such theorems, which we list here to make our article self-contained.

Let $B_r(x)$ be the Euclidean ball with radius $r$ centered at $x$. Then one has the following
\begin{theorem} \label{lsbg} (Lebesgue Differentiation Theorem, see \cite{Lebesgue} and \cite{Rudin:1987:RCA:26851})\\ Let $f:\mathbb{R}^n \rightarrow \mathbb{R}$ be integrable. Then the following is true almost everywhere:
	\[
	f(x) = \lim_{r \rightarrow 0} \frac{1}{\mu (B_r(x))} \int_{B_r(x)} f d\mu,
	\]
	where $\mu$ is the Lebesgue measure on $\mathbb{R}^n$. The above equality holds everywhere if $f$ is continuous on  $\mathbb{R}^n$.
\end{theorem}

Note that this result holds even if we consider balls coming from another equivalent metric structure on $\mathbb{R}^n$, for example if  $B_r(x)$  is an $n$-dimensional cube centered at $x$.  In fact,  the family of balls described above can be replaced by a fairly large family of open sets that {\em ``shrink to $x$ nicely'', } as explained in  \cite{Rudin:1987:RCA:26851}.

Before stating the second theorem, let us recall a definition from measure theory.
\begin{definition}\label{def-abscon}
Let $\mu$ be a positive measure defined on a $\sigma$-algebra $\mathfrak{M}$, and let $\nu$ be a signed measure on $\mathfrak{M}$. Then $\nu$ is called absolutely continuous with respect to $\mu$, denoted by 
$$
\nu\ll\mu,
$$
if $\nu(E)=0$ whenever $\mu(E)=0$.
\end{definition}

\begin{theorem} \label{Rad-Nik} (Radon-Nikodym Theorem, see \cite{Rudin:1987:RCA:26851})\\
If $\mu$ is the Lebesgue measure on $\mathbb{R}^n$ and $\nu$ is a signed measure on $\mathbb{R}^n$ such that $\nu \ll \mu$, then there is a unique integrable real valued function f on $\mathbb{R}^n$ such that for every measurable set A,
\[
 \nu(A) = \int_A f \, d\mu.
\]
$f$ is called the Radon-Nikodym derivative of $\nu$ with respect to $\mu$.

Furthermore, if $\nu$ is a measure (nonnegative) then the function $f$ will be a nonnegative function.
\end{theorem}

For more details about the ``nicely shrinking sets'', Lebesgue differentiation and their relation to the Radon-Nikodym Theorem we refer the reader to Chapter 7 of \cite{Rudin:1987:RCA:26851}.


\subsection{Cone Differentiation Theorem}
In this subsection we derive another generalization of the FTC, which we use later to study the properties of the VLT and its generalizations.

We start with the 2D case. Assume $f$ is an integrable function with respect to the Lebesgue measure on $\mathbb{R}^2$. Let $F(x)$ be defined as in formula (\ref{FTC-n}) with respect to the partial order corresponding to the positive cone generated by some fixed vectors $u,v$.

Define $A_{t,s}(x)$ as the average of $f$ over the parallelogram $P$ centered at $x$, with sides of length $t,s$ and directions $u,v$. We can consider these parallelograms as a family of nicely shrinking  neighborhoods described in \cite{Rudin:1987:RCA:26851}.  Note that the area of the parallelogram made with vectors $tu, sv$ is equal to $\left| \det{(tu,sv)} \right| = ts \left| \det(u,v) \right|$ and we have
\[
A_{s,t}(x) = \frac{1}{ts \left| \det(u,v) \right|} \int_{P} f\, d\mu.
\]

\begin{figure}[h]
\begin{center}
\includegraphics[width=35mm,keepaspectratio]{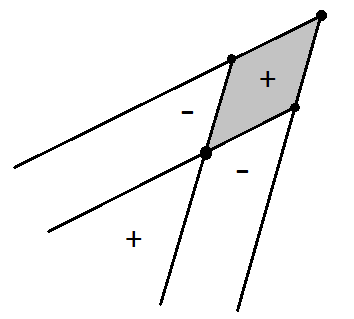}
\end{center}
\caption{Combining the integrals of $f$ over negative cones to get its integral over a parallelogram.}
\label{fig: Vts}
\end{figure}

 Using a simple geometric argument (see Figure \ref{fig: Vts}) and the fact that $F(x)$ is the integral of $f$ over the negative cone at $x$ we get

\begin{multline*}
 A_{t,s} (x) =  \frac{1}{ts \left| \det(u,v) \right|} \Bigg[F\left(x+\frac{t}{2}u+\frac{s}{2}v\right) - F\left(x-\frac{t}{2}u+\frac{s}{2}v\right)  \\
 - F \left(x+\frac{t}{2}u-\frac{s}{2}v\right)	+ F \left(x-\frac{t}{2}u-\frac{s}{2}v\right) \Bigg].
\end{multline*}

Likewise for the $n$-dimensional case, using a geometric argument and induction over $n$, we get the following averaging formula for $f$:

\begin{equation}\label{Vtn}
A_{t_1,\dots,t_n}(x) =
\frac{\sum\limits_{\sigma \in \left\{-\frac{1}{2},\frac{1}{2} \right\}^n } \sgn( \sigma_1 \dots \sigma_n)\, F(x+\sigma_1 t_1 v_1+\dots+ \sigma_n t_n v_n )}
{t_1\dots t_n\, \left| \det(v_1,\dots,v_n) \right|}.
 \end{equation}

In the special case when $t_1 = \dots = t_n = t$, this quantity corresponds to the average of $f$ over $P_t$, the parallelepiped with sides of length $t$ centered at $x$, and we denote it by $A_t(x)$, i.e.
\begin{equation}\label{eq-At}
A_{t}(x) = \frac{1}{\mu(P_t)} \int_{P_t} f\, d\mu =\frac{1}{t^n \left| \det(v_1,\dots,v_n) \right|} \int_{P_t} f\, d\mu.
\end{equation}
Averaging over such infinitesimal symmetric neighborhoods of $x$ and applying Theorem \ref{lsbg} we obtain the following result:

\begin{theorem} \label{ftc}
	Let $\leq$ be an order structure in $\mathbb{R}^n$ corresponding to the positive cone generated by vectors $v_1, \dots ,v_n$. Let $f \in L^1(\mathbb{R}^n)$ and $F(x)$ be defined as in formula (\ref{FTC-n}).
	Then for almost every $x$ we have
\begin{equation}\label{eq_Vtt}
	   f(x) = \lim\limits_{t \rightarrow 0}  A_{t}(x),
\end{equation}
	where $ A_{t}(x)$ is defined in formula (\ref{eq-At}) and can be computed using $t_1=\ldots=t_n=t$ in formula (\ref{Vtn}).
\end{theorem}

Note that this method of recovering $f$ from $F$ is of practical significance, because it is both efficient and simple.\\


Before introducing the main theorem of this section, we  prove one more technical result. In essence, it is a special case of Fubini's theorem, but it is used multiple times throughout the paper, so we prove it here and give it a geometrically descriptive name.

\begin{figure}[h]
\begin{center}
\includegraphics[width=70mm,keepaspectratio]{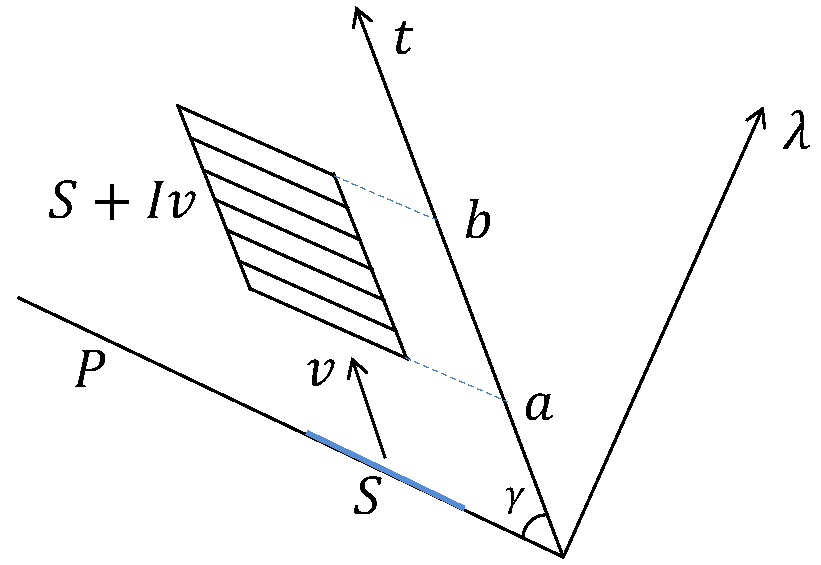}
\end{center}
\caption{A sketch of the sections of $S+Iv$.}
\label{fig: mov_sec}
\end{figure}

Let $P$ be a hyperplane in $\mathbb{R}^n$ and $S\subset P$ be a measurable set. For a vector $v\in\mathbb{R}^n$ transversal to $P$ we denote $S+tv=\{s+tv: s\in S\}$, and for $I=[a,b]\subset \mathbb{R}$ we let $S+Iv=\{S+tv:\, t \in [a,b]\}$. Geometrically, $S+tv$ is an $n-1$ dimensional section of the set $S+Iv$ (see Figure \ref{fig: mov_sec}).

\begin{lemma}\label{mslemma} (Moving Sections Lemma)
Let $f\in L^1(\mathbb{R}^n)$ and $P,S,v$ be defined as above. Then
\[
\int_{S+Iv} f \, d\mu = \sin\gamma  \int_a^b \int_{S+tv} f\, d\mu_P\,  \,dt
\]
where $\mu_P$ is induced from the natural measure on $P$ and $\gamma$ is the angle between $v$ and $S$.
\end{lemma}

\noindent \textbf{Proof.} Let $\lambda$ be the coordinate measured along the axis normal to $P$. Then the Lebesgue measure on $\mathbb{R}^n$ is the product of the Lebesgue measure on the $\lambda$-axis and the natural Lebesgue measure on $P$. Now let $t$ be the coordinate measured along the axis in $v$-direction. Then we have $d\lambda=\sin\gamma\,dt$.

By Fubini's theorem
\[
\int_{S+Iv} f \, d\mu = \int_{a\sin\gamma}^{b\sin\gamma} \int_{S+\lambda \sin\gamma\, v} f\, d\mu_P\, d\lambda  = \int_a^b \int_{S+tv} f\, d\mu_P\, \sin\gamma \,dt.\;\;\;\;\;\; \blacksquare
\] 

In particular, as a consequence of this lemma in $\mathbb{R}^2$ we have the following. Let $S$ be a line segment in $\mathbb{R}^2$ and $v$ be a vector transversal to $S$. If $h(t)=\int_{S+tv} f\, dl$ is the integral along the section $S+tv$, then $\sin\gamma \int_a^b h(t)\,dt$ is equal to the integral of $f$ over the parallelogram $S+Iv$, i.e.
\[
 \int_{S+Iv} f\, d\mu = \sin\gamma \int_a^b h(t)\,dt.
\]

Let us now state the main result of this section.


\begin{theorem} \label{ftccon} (Cone Differentiation Theorem)
Let $\leq$ be an order structure in $\mathbb{R}^n$ corresponding to the positive cone generated by unit vectors $v_1, \dots ,v_n$.  Assume $f\in C_c(\mathbb{R}^n)$ and $F(x) = \int_{y\leq x} f(y) d\mu$. Then we have
\begin{equation}\label{Cone_partial_n}
	 f(x) = \frac{1}{\left| \det(v_1,\dots, v_n) \right|}\, \frac{\partial}{\partial v_1}\dots  \frac{\partial}{\partial v_n} F(x),
\end{equation}
where $\displaystyle\frac{\partial}{\partial v_j}$ is the directional derivative in the direction of $v_j$.
\end{theorem}

\noindent \textbf{Proof.} We provide two different proofs here. The first one is a geometrically intuitive argument applicable to functions in $\mathbb{R}^2$, which is the most relevant case for imaging applications described in the introduction. The second proof is more general and covers arbitrary dimensions.

{\em Proof 1:} Assume $f\in C_c(\mathbb{R}^2)$, and let $u$ and $v$ be the unit vectors generating the positive cone in $\mathbb{R}^2$. For $t\in \mathbb{R}$  and $x\in \mathbb{R}^2$ we define
\[
 q_x(t) = \int_{-\infty}^{0} f(x + tu + sv)\, ds.
\]
 In geometric terms, $q_x(t)$ is the integral of $f$ over the ray propagating from the point $x + tu$ in the direction opposite to $v$. Since $f\in C_c(\mathbb{R}^2)$, one can conclude that $q_x$ is also continuous, e.g. by using the dominated convergence theorem with the dominant $L^1$ function $\chi_{\operatorname{supp}f}\max{|f|}$, where $\chi_{\operatorname{supp}f}$ is the characteristic function of the support of $f$.

Recall that $F$ is the integral  of $f$ over the negative cone. Hence, by the Moving Sections Lemma, we have
\[
F(x) = \sin{(2\beta)} \int_{-\infty}^{0} q_x(t)\, dt,
\]
where $2\beta$ is the angle between $u$ and $v$. Now, using the Lebesgue Differentiation Theorem and continuity of $q_x$ we obtain
$$
\frac{\partial F(x)}{\partial u} = \lim_{h \rightarrow 0} \frac{F(x+hu)-F(x)}{h} \\
 = \lim_{h \rightarrow 0} \frac{\sin{(2\beta)}}{h} \int_0^{h} q_x(t)\, dt \\
 = \sin{(2\beta)}\; q_x(0).
$$
Using the fact that $\sin{(2\beta)}=|\det{(u,v)}|$ and applying the Lebesgue Differentiation Theorem once again we get

\[
 \frac{1}{|\det{(u,v)}|}
\frac{\partial^2 F(x)}{\partial v\, \partial u} =  \lim_{h \rightarrow 0} \; \frac{q_{x+hv}(0) - q_x(0)}{h} = \lim_{h \rightarrow 0} \frac{1}{h} \int_0^{h} f(x+tv)\, dt =  f(x).
\] \\

{\em Proof 2:}
Consider the linear transformation $\phi:\mathbb{R}^n \to \mathbb{R}^n$ defined by $\phi(e_i)=v_i$, $i=1,\ldots, n$, where $\{e_i\}_{i=1}^n$ is the standard basis for $\mathbb{R}^n$. In the domain of $\phi$ consider the partial order corresponding to the positive cone generated by $\{e_i\}_{i=1}^n$. In the range of $\phi$ consider the partial order corresponding to the positive cone generated by $\{v_i\}_{i=1}^n$.

Denote by $C(\hat{x})=\{\hat{y}\in\mathbb{R}^n \, | \, \hat{y}\leq \hat{x}\}$ the negative cone at $\hat{x}$ in the domain of $\phi$. Then $\phi(C(\hat{x}))=\{y\in\mathbb{R}^n,\,y\le x\}$ is the negative cone at $x=\phi(\hat{x})$ in the range of $\phi$.

Using the change of variables formula for $n$-dimensional integrals we obtain
\[
F(x)=\int_{y\le x} f(y)\; dy =
\int_{\phi(C(\hat{x}))} f(y)\; dy =
\int_{C(\hat{x})} f(\phi(\hat{y}))\; |\det(\phi)|\; d\hat{y}.
\]
Now for $x_0 = \phi(\hat{x}_0)$ we have
\begin{equation} \label{eq1}
\begin{split}
\frac{\partial^n F(x_0)}{\partial v_1 \dots \partial v_n}  & = \frac{\partial^n }{\partial v_1 \dots \partial v_n}\bigg|_{x_0} \int_{y\le x} f(y)\, dy \\
 & = \frac{\partial^n }{\partial e_1 \dots \partial e_n}\bigg|_{\hat{x}_0} \int_{C(\hat{x})} f(\phi(\hat{y}))\; | \det(\phi)|\; d\hat{y} \\
  & =  |\det(\phi)|\; \frac{\partial^n }{\partial \hat{x}_1 \dots \partial \hat{x}_n}\bigg|_{\hat{x}_0} \int_{-\infty}^{\hat{x}_n} \dots \int_{-\infty}^{\hat{x}_1} f(\phi(\hat{y}))\; d\hat{y}_1\dots d\hat{y}_n  \\
 & = |\det(\phi)| \; f(\phi(\hat{x}_0)) \\
 & = |\det(v_1, \dots, v_n) |\; f(x_0),
\end{split}
\end{equation}
where we used Fubini's theorem and the Fundamental Theorem of Calculus ($n$ times) correspondingly in the third and fourth lines. Notice, that here $\frac{\partial }{\partial \hat{x}_j }$ denotes the partial derivative with respect to the variable $\hat{x}_j$, while $\frac{\partial}{\partial v_j}$ and $\frac{\partial}{\partial e_j}$ denote the directional derivatives in the direction of vectors $v_j$ and $e_j$ correspondingly. $\hfill\blacksquare$

\vspace{1mm}

\begin{corollary}
In $\mathbb{R}^2$, the Cone Differentiation Theorem can be written as:
\begin{equation}\label{Cone_partial_2}
   f(x) = \frac{1}{\left| \det(u,v) \right|}\; \frac{\partial}{\partial u}  \frac{\partial}{\partial v} F(x).
\end{equation}
\end{corollary}

\begin{remark}
Notice, that Theorem \ref{ftccon} does not imply that $F\in C^n{(\mathbb{R}^n)}$. For example, $\frac{\partial^2}{\partial v^2}F$ may not exist in $\mathbb{R}^2$. \end{remark}

\vspace{2mm}


\section{V-Line Radon Transforms}\label{Chap_BRT}

\begin{figure}[h]
\begin{center}
\includegraphics[width=50mm,keepaspectratio]{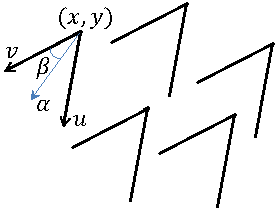}
\end{center}
\vspace{-5mm}
\caption{The $V$-lines in $\mathbb{R}^2$ generated by unit vectors $u,v$.}
\label{fig: V-lines}
\end{figure}

In this section we consider various $V$-line transforms, which map a function on $\mathbb{R}^2$ to its weighted integrals along $V$-shaped trajectories with a fixed axis of symmetry $\alpha$ and a fixed opening angle $2\beta$ (see Figure \ref{fig: V-lines}). We use the \textit{unit} vectors $u, v$ in the directions of the rays of the $V$-lines as generators of the negative cone, which defines a partial order structure of the underlying Euclidean space.

For $(x,y) \in \mathbb{R}^2$ we denote by $R_u(x,y)=\{ (x,y)+tu:\, t \geq 0 \}$ the ray emanating from $(x,y)$ in the direction of $u$. Then the unique V-line with a vertex at $(x,y)$ can be represented by the union $L(x,y)=R_u(x,y) \cup R_v(x,y)$.

In the definitions below we assume that $f\in L^1(\mathbb{R}^2)$ and the transformations are defined for almost every $(x,y)\in\mathbb{R}^2$ (for more details see Appendix II). With some additional regularity assumption on $f$ (e.g. continuity and compact support) the transformations will be defined at every $(x,y)\in\mathbb{R}^2$.

\begin{definition}
 The weighted $V$-line transform $T_{w}$ of $f$ is defined by:
 \begin{equation}\label{def_WV-line}
   (T_wf)(x,y) = c_v \int_{R_v(x,y)} f\, dl + c_u \int_{R_u(x,y)} f\, dl,
 \end{equation}
 where $c_u \ne 0$ and $c_v>0$ are some constants, and $dl$ is the standard Lebesgue measure on the line.
\end{definition}

We introduce additional notations for two special cases, $c_u=c_v=1$ and $c_v=-c_u=1$ .

\begin{definition}
 The (ordinary) $V$-line transform $T$ of $f$ is defined by:
 \begin{equation}\label{def_V-line}
   (Tf)(x,y) = \int_{L(x,y)} f\, dl.
 \end{equation}
 \end{definition}

 \begin{definition}
 The signed $V$-line transform $T_s$ of $f$ is defined by:
 \begin{equation}\label{def_SV-line}
   (T_s f)(x,y) = \int_{R_v(x,y)} f\, dl - \int_{R_u(x,y)} f\, dl.
 \end{equation}
 \end{definition}

In the following subsections we present inversion formulae for these transforms and describe some important properties. We do so by using the techniques of cone differentiation and integration developed in the previous section and some other simple geometric ideas.


\subsection{Inversion of the Ordinary VLT}\label{subsec-IVLT}

Using the Moving Sections Lemma, it is easy to prove that having the VLT of $f$ one can generate its integrals over the corresponding negative cones in $\mathbb{R}^2$. More specifically,

\begin{theorem} \label{inversion1} Let $F(x,y)$ be the integral of $f\in L^1({\mathbb{R}^2)}$ over the negative cone at $(x,y)$. Then
\begin{equation}\label{F_from_T}
  F(x,y) = \sin \beta\int_{0}^{\infty} (Tf)( x+t \alpha_x , y+t \alpha_y ) \, dt.
\end{equation}
\end{theorem}

Using the Cone Differentiation Theorem we immediately obtain an inversion formula for the VLT. Namely, since $|\det(u,v)|=\sin(2\beta)$ we get

\begin{corollary} Let $f\in C_c({\mathbb{R}^2)}$. Then
\begin{equation}\label{VRT-inversion}
 f(x,y) = \frac{1}{2\cos\beta}\, \frac{\partial}{\partial u}  \frac{\partial}{\partial v} \int_{0}^{\infty} (Tf)( x+t \alpha_x , y+t \alpha_y )\, dt.
\end{equation}
If it is only known that $f\in L^1(\mathbb{R}^2)$, then we have the following inversion formula for a.e. $(x,y)$:
\begin{equation}\label{VRT-inversion-2}
f(x,y)  = \lim\limits_{t \rightarrow 0}  A_{t}(x,y),
\end{equation}
where $A_t(x,y)$ is defined as in Theorem \ref{ftc}.
\end{corollary}

Notice that in formula (\ref{VRT-inversion}) the directional derivatives are taken in the direction of the generators of the negative cone, while the Cone Differentiation Theorem uses the generators of the positive cone. But since we have a pair of such derivatives, the algebraic signs appearing due to this difference cancel each other.

\begin{remark}
Formulae that are similar, or equivalent to (\ref{VRT-inversion}), have been obtained in \cite{Florescu-Markel-Schotland, Gouia_Amb_V-line, Kats_Krylov-13, Sherson} using different techniques than those presented herein.
\end{remark}

\begin{figure}[h]
\begin{center}
\includegraphics{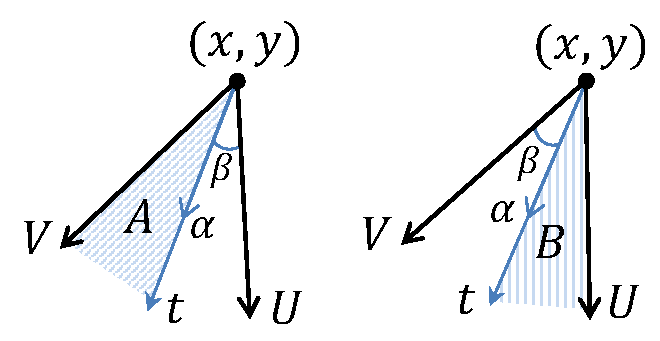}
\end{center}
\caption{The regions $A$ and $B$ with moving sections.}
\label{fig:regions}
\end{figure}

\noindent \textbf{Proof (of Theorem \ref{inversion1}).}
Let $g=Tf$ and $\gamma(t)$ be the parametric equation of the ray starting at $(x,y)$ and moving in the direction of $\alpha=(\alpha_x, \alpha_y)$, i.e. $\gamma(t) = (x,y)+t \alpha$. This ray divides the region enclosed by the $V$-line into two parts: $A$ and $B$ (see Fig. \ref{fig:regions}).

We apply the Moving Sections Lemma to get the integral of $f$ over these regions. In particular,
if we use as moving sections $A_t=V+t\alpha$, then we obtain
\[
\int_A f d\mu = \sin\beta  \int_0^{\infty} \int_{A_t} f\, d\mu_V\, \, dt.
\]
Similarly, for $B_t=U+t\alpha$ we get
\[
\int_B f d\mu = \sin\beta \int_0^{\infty} \int_{B_t} f\, d\mu_U\, \, dt.
\]
Adding these two equations we get the formula
\[
F(x,y) = \int_{A \cup B} f d\mu = \sin\beta \int_0^{\infty}\left( \int_{A_t} f d\mu_U + \int_{B_t} f\, d\mu_V \right)\,  \, dt
\]
\[
 = \sin\beta \int_0^{\infty} g(x+t\alpha_x,y+t \alpha_y)\, \, dt.
\] $\hfill\blacksquare$ \\



\subsection{Inversion of the Weighted VLT}
\begin{figure}[h]
\begin{center}
\includegraphics{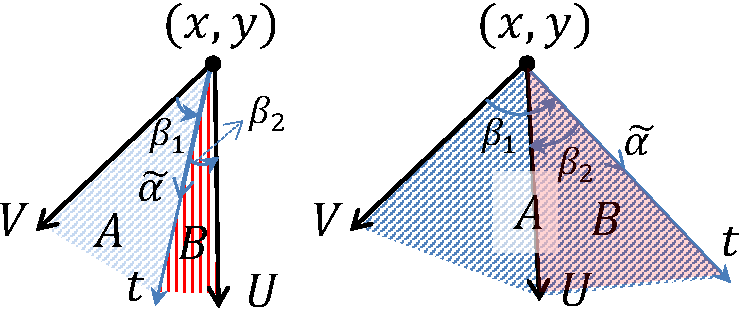}
\end{center}
\caption{The new direction of integration $\widetilde{\alpha}$ and the moving sections. The sketch on the left depicts the setup when $c_u>0$, while the one on the right corresponds to $c_u<0$.}
\label{fig: wregion}
\end{figure}

Consider the weighted V-line transform $T_w$ defined by equation (\ref{def_WV-line}) with constants $c_u\ne 0$ and $c_v>0$. Let us express the fixed opening angle $2\beta$ of the V-line of integration as a sum of two directed angles: $2\beta=\beta_1+\beta_2<\pi$ so that
\begin{equation}\label{def_beta_12}
    \frac{\sin\beta_1}{\sin\beta_2} = \frac{c_v}{c_u}
\end{equation}
(see Figure \ref{fig: wregion}). It is easy to notice that the above relation uniquely defines the angles $\beta_1\in(0,\pi)$ and $\beta_2$ so that $|\beta_2|\in(0,\pi)$. E.g. one can use the following identity
\[
 \cot\beta_1 = \frac{c_u}{c_v \sin(2\beta)} + \cot(2\beta),
\]
and the fact that $\cot x$ is one-to-one on $(0,\pi)$.

Now let $\widetilde{\alpha}=(\widetilde{\alpha}_x, \widetilde{\alpha}_y)$ be the unique unit vector starting from the vertex of the V-line and satisfying the properties:
\begin{equation}
\begin{split}
\operatorname{angle}{(v, \widetilde{\alpha})} = \beta_1, \\
\operatorname{angle}{(\widetilde{\alpha},u)} = \beta_2.
\end{split}
\end{equation}

The vector $\widetilde{\alpha}$ can also be expressed as a linear combination of the unit vectors $u$ and $v$ as follows:
\begin{equation}\label{alpha_tilde}
\widetilde{\alpha} = \frac{c_u v + c_v u}{\|c_u v + c_v u\|}.
\end{equation}
To verify the last relation we use the cross products:

\[
 \sin {\beta_1} = \| v \times \widetilde{\alpha} \| =  \frac{\| v \times (c_u v + c_v u) \|}{\|c_u v + c_v u\|} =
 \frac{c_v \sin{(2\beta)}}{\|c_u v + c_v u\|},
\]

\[
 |\sin {\beta_2}| = \| \widetilde{\alpha} \times u \| =  \frac{\| (c_u v + c_v u) \times u \|}{\|c_u v + c_v u\|} =
 \frac{|c_u| \sin{(2\beta)}}{\|c_u v + c_v u\|}.
\]
Notice, that formula (\ref{alpha_tilde}) implies that $\operatorname{sgn}{(\sin \beta_2)}=\operatorname{sgn} {(c_u)}$. As a result
\[
 \sin {\beta_2} =  \frac{c_u \sin{(2\beta)}}{\|c_u v + c_v u\|},
\]
and formula (\ref{def_beta_12}) follows immediately.

Just as in the case of the ordinary VLT, to invert the weighted VLT we express $F(x,y)$ through $T_wf$ and apply the Cone Differentiation Theorem.

\begin{theorem}\label{arb_w_thm}
Let $F(x,y)$ be the integral of $f\in L^1({\mathbb{R}^2)}$ over the negative cone at $(x,y)$. If $c_v,\, c_u,\, \beta_1,\, \beta_2$, and $\widetilde{\alpha}$ are defined as above, then
\begin{equation}\label{F_from_Tw}
F(x,y) = \frac{\sin\beta_1}{c_v} \int_{0}^{\infty} (T_wf)(x+\widetilde{\alpha}_x t,y+\widetilde{\alpha}_y t)\, dt.
\end{equation}
\end{theorem}
\noindent \textbf{Proof.}
We follow the same steps as in the proof of the non-weighted version, but this time we divide the integration region using a new direction $\widetilde{\alpha}$ defined in the statement of the theorem.

Applying the Moving Sections Lemma we obtain:
\[
\int_A f\, d\mu = \sin\beta_1  \int_{0}^{\infty} \int_{V+t \widetilde{\alpha}} f \,d\mu_V \, dt,
\]
\[
\int_B f\, d\mu = \sin\beta_2 \int_{0}^{\infty} \int_{U+t\widetilde{\alpha}} f \,d\mu_U  \, dt.
\]

Notice, that if $c_u>0$ then $F(x,y) = \int_{A \cup B} f\, d\mu$, while in the case of $c_u<0$ we have $F(x,y) = \int_{A \setminus B} f\, d\mu$.

In both cases we get
\[
F(x,y) =  \int_{0}^{\infty} \left( \sin\beta_1\int_{V+t\widetilde{\alpha}} f\, d\mu_V + \sin\beta_2 \int_{U+t \widetilde{\alpha}} f \, d\mu_U \right) \, dt.
\]
Using the relation $c_v/c_u = \sin\beta_1/\sin\beta_2$ we have
\[
F(x,y)=\frac{\sin\beta_1}{c_v} \int_{0}^{\infty} \left( c_v\int_{V+t\widetilde{\alpha}} f\, d\mu_V + c_u \int_{U+t\widetilde{\alpha}} f \, d\mu_U \right) \, dt
\]
\[
=\frac{\sin\beta_1}{c_v} \int_{0}^{\infty} (T_wf)(x+\widetilde{\alpha}_x t, y+\widetilde{\alpha}_y t)\, dt,
\]
which finishes the proof.
$\hfill\blacksquare$ \\

In the special case of the signed VLT (i.e. $c_v=-c_u=1$), we have $\widetilde{\alpha} = \frac{u-v}{\| u-v\|}$, which is the unit vector perpendicular to the cone direction $\alpha$. As a result we get the following

\begin{theorem}
Let $F(x,y)$ be the integral of $f\in L^1({\mathbb{R}^2)}$ over the negative cone at $(x,y)$ and let $T_s f$ be the signed V-line transform of $f$. For $\widetilde{\alpha} = \frac{u-v}{\| u-v\|}$ we have
\begin{equation}
F(x,y) = \cos \beta \int_{0}^{\infty} (T_sf)(x+\widetilde{\alpha}_x t,y+\widetilde{\alpha}_y t)\, dt.
\end{equation}
\end{theorem}

Applying the Cone Differentiation theorem to Theorem \ref{arb_w_thm} one can now get an inversion formula for the weighted VLT.

\begin{theorem}
Let $T_w f$ be the weighted V-line transform of $f\in C_c({\mathbb{R}^2)}$, with arbitrary non-zero weights. For $\widetilde{\alpha} =\displaystyle \frac{c_u v + c_v u}{\|c_u v + c_v u\|} $ we have
\begin{equation}\label{WVLT-inv}
 f(x,y) = \frac{1}{\|c_u v + c_v u\|}\, \frac{\partial}{\partial u}  \frac{\partial}{\partial v} \int_{0}^{\infty} (T_wf)(x+\widetilde{\alpha}_x t,y+\widetilde{\alpha}_y t)\, dt.
\end{equation}
The coefficient in the above formula can be expressed as
\begin{equation}\label{nrm_alph}
\frac{1}{\|c_u v + c_vu\|}=
\displaystyle\frac{\sin{\beta_1}}{c_v\sin{(2\beta)}}.
\end{equation}
\end{theorem}

\vspace{5mm}

Assume that $f\in C_c(\mathbb{R}^2)$ and its values are known on the boundary of some bounded, open, convex set $\Omega$. Then one can recover $f$ inside $\Omega$ using its weighted VLT data restricted to V-lines with vertices inside $\Omega$. Namely

\begin{theorem}\label{bound_inv}
Consider a bounded, open, convex set $\Omega$ in $\mathbb{R}^2$ and let $g=T_w f$ be the weighted V-line transform of $f\in C_c(\mathbb{R}^2)$. For each point $p=(x,y)\in \Omega$ let
\[
m_p=\inf\,\{s\ge 0 : p+s\widetilde{\alpha} \in \partial\Omega \},
\]
\[
m_n = m_p-1/n,
\]
\[
(x_n,y_n) = (x,y) + m_n \widetilde{\alpha},
\]
\[
(x_0,y_0) = (x,y) + m_p \widetilde{\alpha}.
\]
Then $(x_n,y_n) \rightarrow (x_0,y_0)$ and
\begin{equation}\label{VRT-inversion2}
f(x,y) = f(x_0,y_0) + \lim\limits_{n\rightarrow \infty} \frac{\sin{\beta_1}}{c_v \sin{(2\beta)}}\, \frac{\partial}{\partial u}\,\frac{\partial}{\partial v}\, \int_0^{m_n}   g(p+s \widetilde{\alpha})\, ds.
\end{equation}
\end{theorem}
\noindent \textbf{Proof.} For large values of $n$ we have $(x_n,y_n)\in \Omega$, since $\Omega$ is convex and the ray 
$\{ p+s\widetilde{\alpha}: s\ge 0 \}$ is transverse to $\partial\Omega$. Hence one can choose a number $\delta(n)>0$ so that the balls $B_1, B_2$ of radius $\delta$ centered at $(x,y)$ and $(x_n,y_n)$ respectively are contained in $\Omega$, i.e. $B_1 \subseteq \Omega, B_2 \subseteq \Omega$. By convexity, any convex combination of points in $B_1$ and $B_2$ belongs to $\Omega$.
By the previous theorem, for any pair $p \in B_1$ and $p'=p+t \widetilde{\alpha} \in B_2$ we know that

\[
\frac{\sin{\beta_1}}{c_v \sin{(2\beta)}}\, \frac{\partial}{\partial u}\,\frac{\partial}{\partial v}\, \int_0^{\infty}   g(p+s \widetilde{\alpha})\, ds = f(p),
\]

\[
\frac{\sin{\beta_1}}{c_v \sin{(2\beta)}}\, \frac{\partial}{\partial u}\,\frac{\partial}{\partial v}\, \int_t^{\infty}   g(p+s \widetilde{\alpha})\, ds = f(p').
\]

Hence,  using the equality
\[
\int_0^{m_n}   g(p+s \widetilde{\alpha})\, ds = \int_0^{\infty}   g(p+s \widetilde{\alpha})\, ds - \int_{m_n}^{\infty}   g(p+s \widetilde{\alpha})\, ds
\]
We get

\[
\frac{\sin{\beta_1}}{c_v \sin{(2\beta)}}\, \frac{\partial}{\partial u}\,\frac{\partial}{\partial v}\, \int_0^{m_n}   g(p+s \widetilde{\alpha})\, ds =  f(x,y) - f(x_n,y_n).
\]
Taking the limit of this equality when $n\to\infty$ completes the proof. $\hfill\blacksquare$
\begin{remark}
In the case of the signed VLT (i.e. $c_v=-c_u=1$) a similar inversion formula was obtained in \cite{Kats_Krylov-13} using other techniques.
\end{remark}

\begin{remark}
In the special case when $c_u=0$, assuming $c_v=1$, the inversion of VLT reduces to a trivial application of the FTC. Notice, that the weighted VLT in this setup is essentially (up to a constant multiple) the same as the ordinary VLT that uses V-lines with an opening angle $2\beta=0$, i.e. with coinciding branches.
\end{remark}


\subsection{A Range Description for the VLT}

We start this subsection with providing the necessary and sufficient conditions for a function $F$ to be a cone integral of another function $f\geq 0 $ with respect to a given order structure in $\mathbb{R}^n$. In other words we answer the question: for which $F$ is there an $f\ge0$ such that $F(x)=\int_{y\leq x} f(y)\,d\mu$? \\

 Our approach is motivated by the corresponding result for $\mathbb{R}^1$ stated in Remark \ref{abs_cont_1}. We use the Radon-Nikodym Theorem to get the desired description of $F$. For an appropriate $F$, we construct a corresponding measure $\nu$, for which Theorem \ref{Rad-Nik} implies the existence of its Radon-Nikodym derivative $f$.

 For $x\in \mathbb{R}^n, c_i \in \mathbb{R}^+$, let $P(x,c_1,\dots ,c_n)$ be an $n$-dimensional half-open parallelepiped defined by
 \begin{equation}\label{parallele}
  P(x,c_1,\dots, c_n) = \left\{ x+ \sum_{i=1}^{n} t_i v_i \in \mathbb{R}^n ; -c_i \leq t_i < c_i \right\},
 \end{equation}
where $v_1, \ldots, v_n$ are the basis vectors defining the order structure in $\mathbb{R}^n$. These parallelepipeds  are the analogs of intervals in $\mathbb{R}^1$.

For a given function $F:\mathbb{R}^n\rightarrow \mathbb{R}$ we define a set function $\nu_0$ on the ring of subsets generated by these parallelepipeds by:
 \[
  \nu_0(P(x,c_1,\dots, c_n)) = \sum_{\sigma \in \{-1,1\}^n } \sgn(\sigma_1 \dots \sigma_n) F(x+\sigma_1 c_1 v_1+\dots + \sigma_n c_n v_n )
 \]
 and extending $\nu_0$ to the ring using the outer measure induced by $\nu_0$.

Using an analogy with absolute continuity on the real line we define absolute continuity of $F$ as follows:
\begin{definition}\label{abscon}
	Let $F$ be a  function on $\mathbb{R}^n$ with a given order structure. Then we say $F$ is {\em absolutely continuous} if and only if for every $\epsilon > 0$ there exists $\delta > 0$ such that for any finite collection of disjoint parallelepipeds $\{P_i \}_{i=1}^m$,   $\sum_{i=1}^{m} \mu(P_i) < \delta$ implies $\sum_{i=1}^{m} | \nu_0(P_i) | < \epsilon$.
\end{definition}

Let us call $\alpha = \displaystyle\frac{v_1+\dots+v_n}{\| v_1+\dots+v_n \|}$ the direction of the negative cone.

\begin{definition}\label{P-increase}
We say that a  function $F$ on $\mathbb{R}^n$ is {\em  P-cumulative } with respect to the given order structure, if

\begin{itemize}
    \item $\nu_0(P) \geq 0$ for every parallelepiped $P$ (non-decreasing condition);
    \item $\lim\limits_{t\rightarrow \infty} F(t\alpha) = 0$ and $\lim\limits_{t\rightarrow -\infty} F(t\alpha) < \infty$.

\end{itemize}

\end{definition}

The second condition may look backwards for readers familiar with cumulative probability distributions. It is simply due to the fact that we use integrals over negative cones and $\alpha$ points in the negative direction.

\begin{remark}
Note, that if $f\geq 0$ is integrable then $F$ defined by formula (\ref{FTC-n}) will be P-cumulative with respect to the underlying order structure.
\end{remark}

If $F$ is P-cumulative, then $\nu_0$ is a pre-measure on the ring of subsets generated by the parallelepipeds. By applying the Caratheodory Extension Theorem, we can extend $\nu_0$ to a measure $\nu$ on $\mathbb{R}^n$ (the domain of $\nu$ is the $\sigma$-algebra of Lebesgue measurable sets). The details of this construction can be found in many measure theory texts (see for example \cite{Athreya}).

We observe that for $F(x) = \int_{t\leq x}f(t)\,d\mu$, if we send $x$ to infinity in the direction opposite to $\alpha$, then in the limit we get the integral of $f$ over $\mathbb{R}^n$. In other words,
\[
 \lim\limits_{t\rightarrow -\infty} F(t\alpha) = \int_{\mathbb{R}^n} f\,d\mu.
\]
   Hence, for nonnegative $f$, if $\lim\limits_{t\rightarrow -\infty} F(t\alpha)  < \infty$ we have $f \in L^1(\mathbb{R}^n)$.

   The following is a consequence of the Radon-Nikodym Theorem:

\begin{theorem}\label{fexist}
	Let $F$ be a P-cumulative function on $\mathbb{R}^n$. Then $F$ is absolutely continuous if and only if there exists a nonnegative function $f$ such that
	\[
	 F(x) = \int_{y\leq x} f \, d\mu.
	\]
\end{theorem}
\noindent \textbf{Proof.} Assume $f\ge0$. Then (as we discussed above) $f$ is integrable and therefore locally integrable. This implies that $f$ is the density of some absolutely continuous measure  $\lambda \ll \mu$, i.e.
$$
\lambda(E)=\int_E f\, d\mu.
$$
Notice, that since $f\in L^1{(\mathbb{R}^2)}$, it follows that $\lambda$ is a finite measure, hence one can use the ``$\epsilon-\delta$'' definition of absolute continuity (see Appendix I). Hence, for any $\epsilon > 0$ there exists a $\delta>0$ such that $\mu(E)<\delta$ implies $\lambda{(E)}<\epsilon$.

Combining the above relation with
\[
  \sum_{i=1}^{m} | \nu_0(P_i) | = \sum_{i=1}^{m} \int_{P_i}f\, d\mu = \int_{\cup P_i}f \, d\mu
\]
we establish that $F$ is absolutely continuous in the sense of Definition \ref{abscon}.

 The other direction of the theorem is an implication of the Radon-Nikodym Theorem. Here we use the fact that the constructed measure $\nu$ is absolutely continuous with respect to $\mu$ if $F$ is absolutely continuous in the sense of Definition \ref{abscon}  (see Appendix I for more details). $\hfill\blacksquare$ \\

 In particular, this shows that if $f$ is related to $F$ as described in the previous theorem, then for any measurable set $S$ we have $\nu(S)=\int_S f\, d\mu$. \\



We will now use the previous theorem to provide a range description for the VLT. Let us start with a new notation and a technical result that will be handy in the proof of the main theorem.

\begin{figure}[h]
\begin{center}
\includegraphics[width=40mm,keepaspectratio]{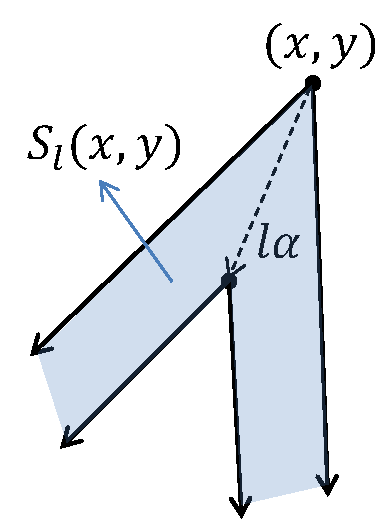}
\end{center}
\caption{The ``broken strip'' $S_l(x,y)$.}
\label{fig:br-strip}
\end{figure}

For $(x,y)\in \mathbb{R}^2$ let $$S_l(x,y)=\{ (x_1,y_1)\, |  (x_1,y_1) \leq (x,y) \} \setminus \{ (x_1,y_1)\, | \, (x_1,y_1) \leq (x,y)+\alpha l \}$$ be the ``broken strip'' region of width $l\sin\beta>0$ (see Figure \ref{fig:br-strip}).

\begin{lemma}\label{lem_bstr-bl}
For $f\in L^1(\mathbb{R}^2)$  and $(x,y)\in \mathbb{R}^2$ we have
\[\lim\limits_{\varepsilon \rightarrow 0}\frac{1}{\varepsilon} \int_{S_{\varepsilon}(x,y)}f\,d\mu = \sin \beta \int_{L(x,y)} f\,dl \]
almost everywhere, where $L(x,y)$ is the V-line at $(x,y)$.
\end{lemma}
\noindent \textbf{Proof.} Define  $k(t) = \int_{L((x,y)+t\alpha)} f \, dl$. By the Moving Sections Lemma we have
\[
\int_{S_{\varepsilon}(x,y)} f \,d\mu = \sin \beta  \int_{0}^{\varepsilon} k(t) \, dt.
\]
Dividing both sides by $\varepsilon$, taking a limit and applying Theorem \ref{lsbg}  (Lebesgue Differentiation Theorem) we get
\[
\lim\limits_{\varepsilon \rightarrow 0}\frac{1}{\varepsilon} \int_{S_{\varepsilon}(x,y)}f\,d\mu =
\sin \beta\;
\lim\limits_{\varepsilon \rightarrow 0}\frac{1}{\varepsilon}
\int_{0}^{\varepsilon} k(t) \, dt=
\sin\beta\; k(0) = \sin \beta \int_{L(x,y)} f\,dl,
\]
for almost every $(x,y)\in\mathbb{R}^2$. The use of Theorem \ref{lsbg} is valid, since $k(t)\in L^1(\mathbb{R})$ by Fubini's theorem applied to $f\in L^1(\mathbb{R}^2)$. $\hfill\blacksquare$ \\

Now we use Theorem \ref{fexist} and the above Lemma to prove the main result of this subsection.

\begin{theorem} (Range Description)
	A function $g$ on $\mathbb{R}^2$
	is the image of some nonnegative function $f\in L^1(\mathbb{R}^2)$ under the V-line transform if and only if the function $F$ defined by
	\[
	 F(x,y) = \sin\beta \int_{0}^{\infty} g( x+t \alpha_x , y+t \alpha_y ) \, dt
	\]
	is P-cumulative and absolutely continuous in the sense of Definition \ref{abscon}.
\end{theorem}
\noindent \textbf{Proof.} Let $g$ be the image of a nonnegative function $f \in L^1(\mathbb{R}^2)$ under the VLT. Then  we have
\[
F(x,y) =\sin \beta \int_{0}^{\infty} \int_{L( x+t \alpha_x , y+t \alpha_y )} f\, dl\,  \, dt = \int_{z\leq (x,y)} f\, d\mu
\]
and hence by Theorem \ref{fexist}, $F$ is absolutely continuous and P-cumulative.\\
For the proof in the other direction, assume that $F$ is absolutely continuous and P-cumulative. By Theorem \ref{fexist} there exists a nonnegative function $f$ such that
\[
F(x,y) = \sin \beta\int_{0}^{\infty} g( x+t \alpha_x , y+t \alpha_y ) \, dt
= \int_{z\leq (x,y)} f\, d\mu.
\]
At the same time Lemma \ref{lem_bstr-bl} implies that for almost every $(x,y)\in\mathbb{R}^2$
\[
 \lim\limits_{\varepsilon \rightarrow 0} \frac{F(x,y)- F(x+\varepsilon \alpha_x, y+\varepsilon \alpha_y)}{\varepsilon} =\sin\beta \int_{L(x,y)} f\, d\mu.
\]
On the other hand, we can rewrite the left hand side of this relation as follows and apply Theorem \ref{lsbg}
\[
 \lim\limits_{\epsilon \rightarrow 0} \frac{1}{\epsilon} \int_{0}^{\epsilon}  g( x+t \alpha_x , y+t \alpha_y ) \sin \beta\, dt = \sin \beta\, g(x,y)
\]
for almost every $(x,y)$. Hence, $g(x,y) = \int_{L(x,y)} f\, d\mu$ almost everywhere. Also, $f$  is integrable because
\[
\int_{\mathbb{R}^2} f \, d \mu = \lim\limits_{t\rightarrow -\infty} F(t\alpha_x, t\alpha_y) =\sin\beta \int_{-\infty}^{+\infty} g(t\alpha_x, t\alpha_y)  \, dt < \infty.
\]
$\hfill\blacksquare$

\begin{theorem} (Weighted Case)
	A function $g$ on $\mathbb{R}^2$
	is the image of a nonnegative function $f\in L^1(\mathbb{R}^2)$ under the weighted VLT if and only if the function $F$, as defined  in Theorem \ref{arb_w_thm},	is P-cumulative and absolutely continuous.
\end{theorem}
\noindent \textbf{Proof.} According to Theorem \ref{arb_w_thm}, $F$ represents the cone integral of $f$ and hence we can apply the same proof as in the previous non-weighted case. $\hfill\blacksquare$

The only other known-to-us range description of the VLT  was given in \cite{Kats_Krylov-13}, but it uses more data (three detectors) and is of an entirely different nature.

\subsection{Support Theorems for the VLT}

\begin{theorem} \label{supp}
	Let $T_wf$ be the weighted V-line transform of $f\in L^1(\mathbb{R}^2)$, and $F$ be defined as in formula (\ref{F_from_Tw}).
	If $F$ is constant on some $S\subset \mathbb{R}^2$ with a non-empty interior $S^0$, then $f=0$ in $S^0$.
\end{theorem}
\noindent \textbf{Proof.} Let $x$ be an interior point of $S$. Then we can find $t$ such that all parallelograms of size smaller that $t$ centered at $x$ will lie inside $S$. Since $F$ has the same values at the corners of these parallelograms, by definition we have $A_{t}(x)=0$. Hence, $f(x) = \lim\limits_{t \rightarrow 0} A_{t}(x) = 0$ by Theorem \ref{ftc}. $\hfill\blacksquare$\\

\begin{theorem} \label{supp2}
	Let $f\in L^1(\mathbb{R}^2)$ be continuous and
	$T_wf=0$ on some $S\subset \mathbb{R}^2$. Let $L$ be a line parallel to the vector $\widetilde\alpha$ defined in (\ref{alpha_tilde}) such that $L\cap S\ne\varnothing$.	Then $f$ is constant on each connected component of $L \cap S$.

In particular, if $S$ is a compact set, $T_wf=0$ in $S$, and $f=0$ on the boundary $\partial S$ of $S$, then $f\equiv0$ in S.
\end{theorem}
\noindent \textbf{Proof.} Take two points on a connected component of $L \cap S$ and let $\Omega$ be the line interval connecting these two points. Then $\Omega$ is a compact convex set and by Theorem \ref{bound_inv} the value of $f$ should be the same at both endpoints.

The second part follows from the fact that for any $x$ in the closed and bounded $S$ the intersection of the ray $L=\{ x+t\alpha | \; t\in \mathbb{R}  \}$ and $\partial S$ is nonempty. $\hfill\blacksquare$\\ 
\section{A Conical Radon Transform in Higher Dimensions}\label{Chap-nd}
Now we consider a generalization of our results to higher dimensions. As before, we assume that $f\in L^1(\mathbb{R}^n)$ and the transformation is defined at almost every $x\in\mathbb{R}^n$.

\begin{figure}[h]
\begin{center}
\includegraphics[width=40mm,keepaspectratio]{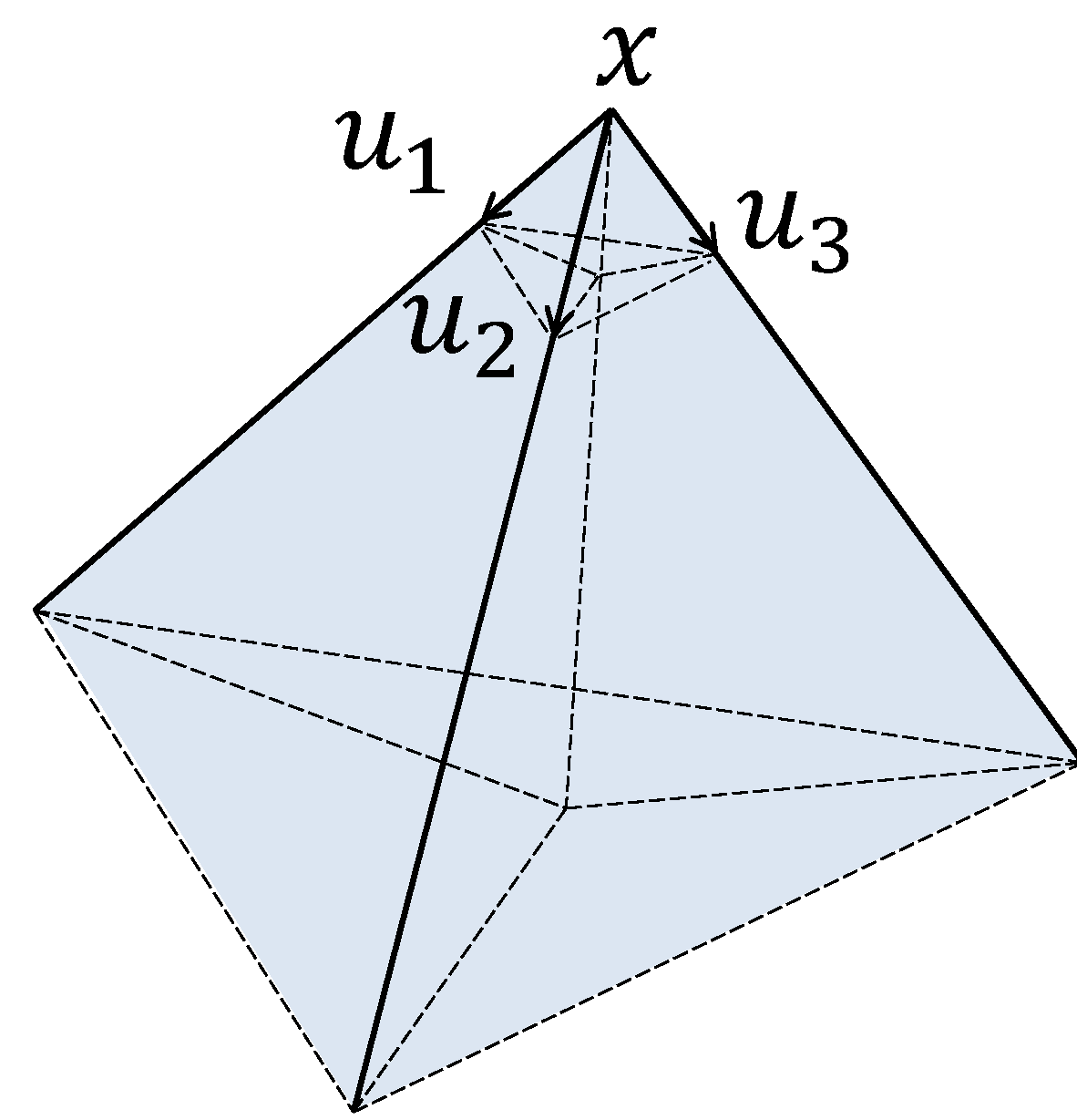}
\end{center}
\caption{A polyhedral cone in $\mathbb{R}^3$ generated by unit vectors $u_1,u_2,u_3$.}
\label{fig:polcone}
\end{figure}

\begin{definition}
The conical Radon transform $T$ maps $f$ into the set of its integrals over the boundaries $\partial C(x)$ of polyhedral cones $C(x)$ generated by fixed unit basis vectors $u_1,\dots,u_n$ starting from $x$ (see Figure \ref{fig:polcone}). Namely, \\
\[
(Tf)(x) = \int_{\partial C(x)} f \,dS,
\]
where $dS$ is the standard $n-1$ dimensional Lebesgue measure on $\partial C$.
\end{definition}

Notice that the number of edges (and faces) of the polyhedral cone coincides with the dimension of the underlying space. In this case there exists a unique unit vector $\alpha$ such that, when starting from the vertex $x$ of the cone, $\alpha$ is pointing inside the cone and has the same angle $\beta$ with all $n$ faces of the cone.

Let $ X_i = \operatorname{span}\langle u_1,\dots,u_{i-1},u_{i+1},\dots,u_n \rangle$ denote the hyperplane containing the $i$-th face of the polyhedral cone and define $y_i$ to be the unit vector in $X_i^{\perp}$ such that
\[
\langle \alpha,y_i \rangle=\sin\beta,\;\;\;  \forall i=1,\dots,n.
\]
The following theorem is the $n$-dimensional analogue of Theorem \ref{inversion1}, as it provides a formula for generating the integral of $f$ over polyhedral cones $C$ from the conical Radon transforms $Tf$.

\begin{theorem}\label{th_con_tr}
	Let $f\in L^1(\mathbb{R}^n)$, and $T,\alpha, y_j$ be defined as above. Then
\begin{equation}\label{def_con_tr}
	F(x) = \langle \alpha,y_1\rangle \int_{0}^{\infty} (Tf)(x+t\alpha)  \, dt
\end{equation}
	is the integral of $f$ over the cone generated by $u_1, \dots, u_n$ with vertex at $x$.
\end{theorem}
\noindent \textbf{Proof.} We use a strategy similar to the two-dimensional case by replacing regions $A,B$ in that proof with $\{A_i\}_{i=1}^{n}$. Let $C(x)$ be the cone at $x$. Then $C = \cup_{i=1}^n A_i$, where $A_i$ is the cone made by $u_1,\dots,u_{i-1},u_{i+1},\dots,u_n $ and $\alpha$. In other words, we break $C$ into disjoint components using the vector $\alpha$.

 Now, to integrate $f$ over $C$ we write
 \[
  \int_C f \, d\mu = \sum_{i=1}^{n} \int_{A_i} f\, d\mu.
 \]
 A straightforward application of the Moving Sections lemma yields:
 $$
 \int_{A_i} f\, d\mu =
 \sin\beta \int_{0}^{\infty}
 \int_{X_i+t\alpha} f \, d\mu_X\;   \, dt.
 $$
 At the same time
 $$
 \sum_{i=1}^{n} \int_{X_i+t\alpha} f\, d\mu_X = Tf(x+t\alpha),
 $$
which finishes the proof. $\hfill\blacksquare$ \\

\begin{corollary}
One can invert the conical Radon transform $T$ by using formula (\ref{def_con_tr}) to generate $F$ from $Tf$ and then applying Theorem \ref{ftc} or Theorem \ref{ftccon}.
\end{corollary}

\begin{remark}
In analogy with the 2D case, once can consider a weighted conical Radon transform, where the integration along each face of the polyhedral cone is done with a different constant weight. An inversion procedure for such a transform can be obtained following the approach of the 2D case and the previous corollary.
\end{remark}

\begin{remark}
Theorems \ref{bound_inv}, \ref{supp}, \ref{supp2} can all be generalized to the case of the conical Radon transform, with proofs following the corresponding arguments used in the case of the V-line transform.
\end{remark}

\begin{remark}
We are not aware of any imaging application for the conical Radon transforms (both circular and polyhedral) with vertices inside the image domain. At this point, the study of such transformations is of purely theoretical interest. The applicability of our approach developed in this paper to the case of circular cones is subject of an ongoing work.
\end{remark} 
\section{Numerical Simulations}\label{Chap_num}

\subsection{Reconstructing a function from its VLT}
We work with a standard 800x800 Shepp-Logan phantom in MATLAB and consider V-lines with an opening angle $\beta = \operatorname{arctan}(1/2)$ and a horizontal symmetry axis. To compute the one-dimensional ray integrals we use a linear interpolation to evaluate the image values along the given line direction with specific step size of dx=0.8 (in pixel). Adding two ray integrals, in two directions, at any given point we get the V-line data.

\begin{figure}[H]
	
	\begin{center}
		
		\begin{subfigure}{.5\textwidth}
			\centering
			\includegraphics[width=.7\linewidth]{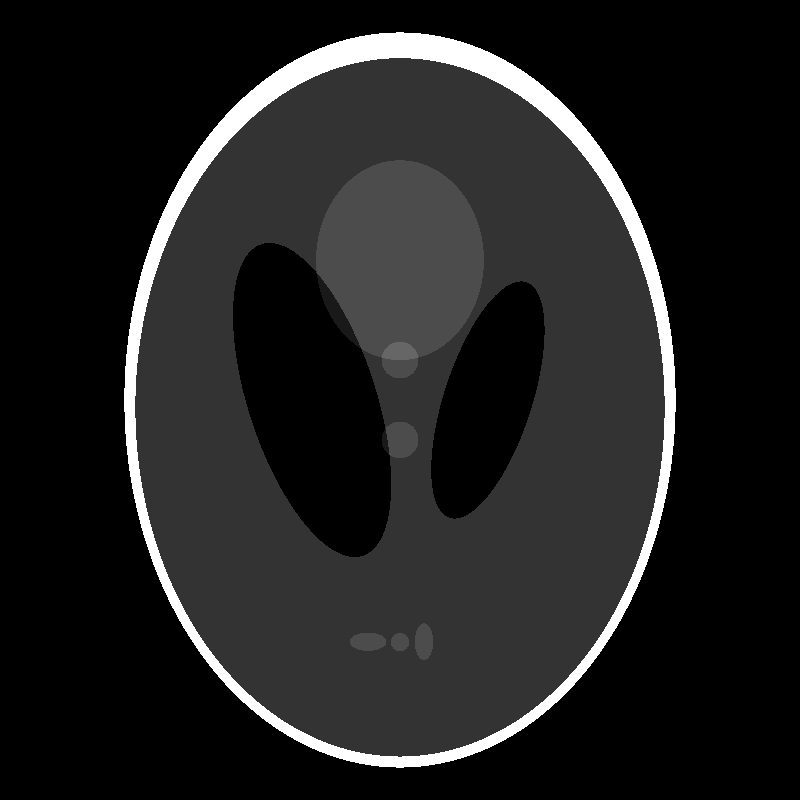}
			\caption{Original Phantom.}
			\label{fig:sfig1}
		\end{subfigure}%
		\begin{subfigure}{.5\textwidth}
			\centering
			\includegraphics[width=.7\linewidth]{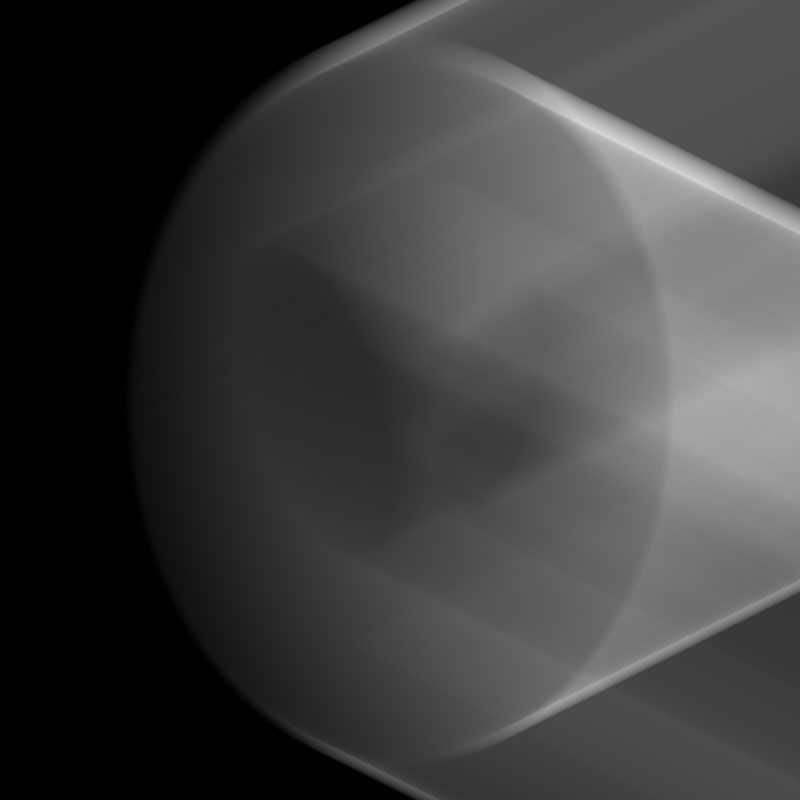}
			\caption{V-Line Data.}
			\label{fig:sfig2}
		\end{subfigure}
		\caption{}
		\label{fig:fig}
		
	\end{center}
	
\end{figure}

We present below two different numerical reconstructions of the phantom from its VLT data. The first one is based on formula (\ref{VRT-inversion}) and uses the directional derivatives of the conical integrals of the image function $f$.
\begin{figure}[H]
	\begin{center}
			\includegraphics[width=.35\linewidth]{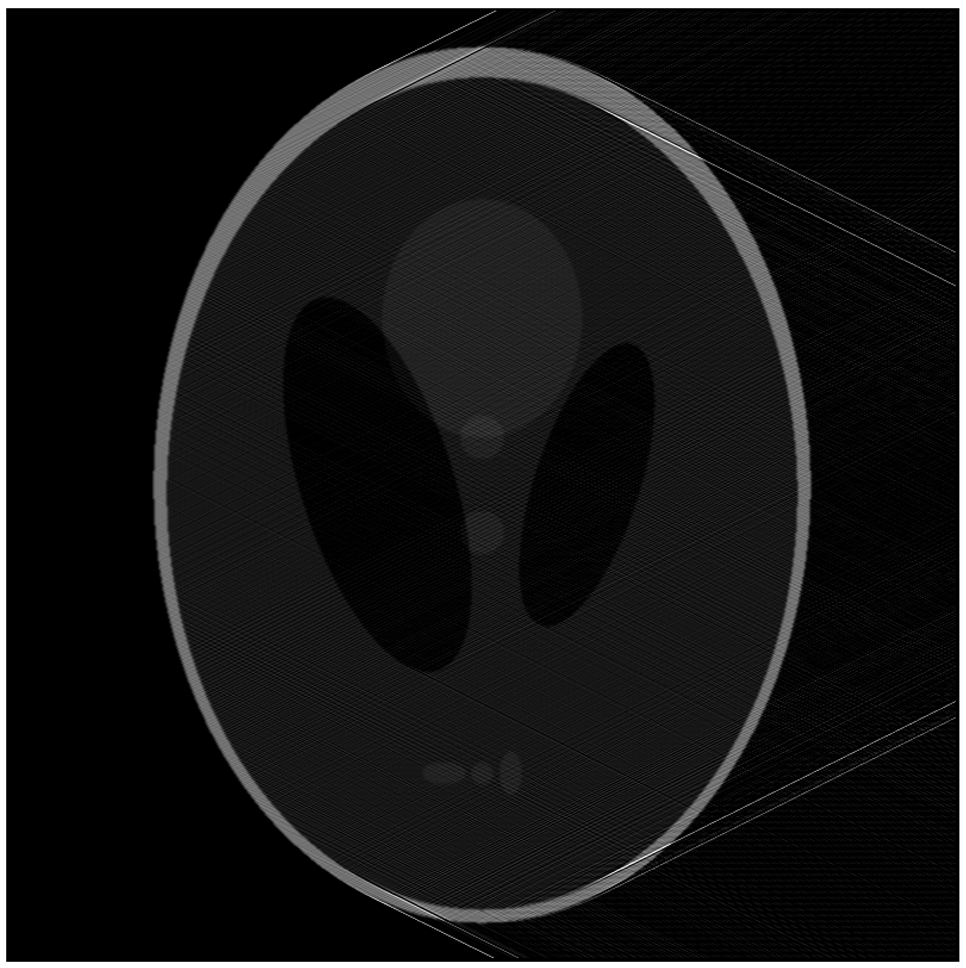}
			\caption{Reconstruction using formula (\ref{VRT-inversion}).}
			\label{fig:CDT}
	\end{center}
\end{figure}

As we have mentioned before, formulae that are very similar or equivalent to (\ref{VRT-inversion}) have been obtained in \cite{Florescu-Markel-Schotland, Gouia_Amb_V-line, Kats_Krylov-13, Sherson} using different techniques than those presented in this paper. Most of those works also included numerical simulations, which agree well with our numerical reconstruction presented in Figure \ref{fig:CDT}.

\begin{figure}[H]
	\begin{center}
	\begin{subfigure}{.5\textwidth}
		\centering
		\includegraphics[width=.7\linewidth]{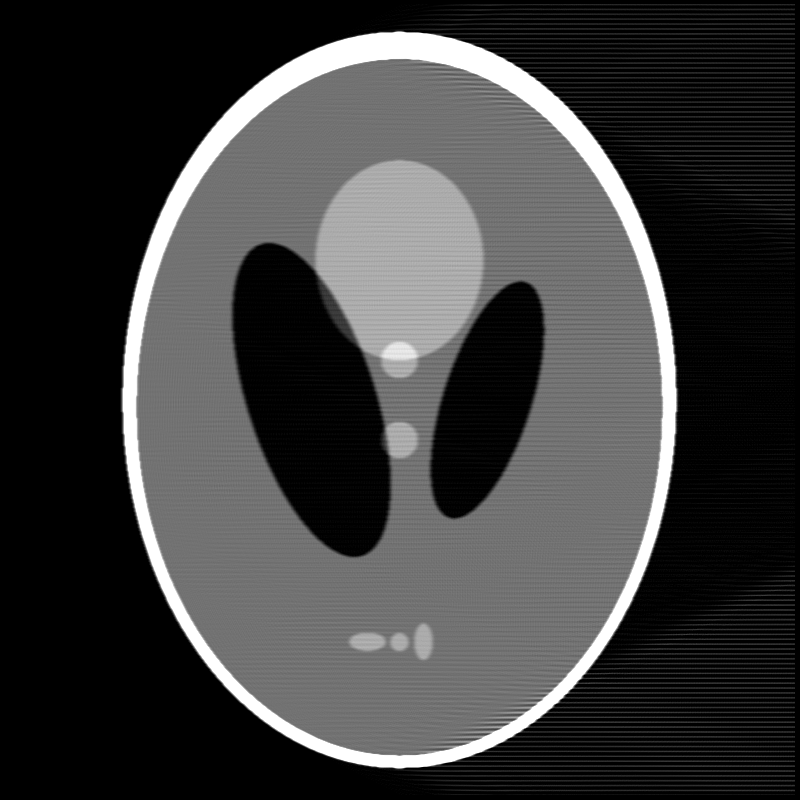}
		\caption{Reconstructed phantom ($\epsilon = 1$).}
		\label{fig:fig2}
	\end{subfigure}%
	\begin{subfigure}{.5\textwidth}
	\centering
		\includegraphics[width=.7\linewidth]{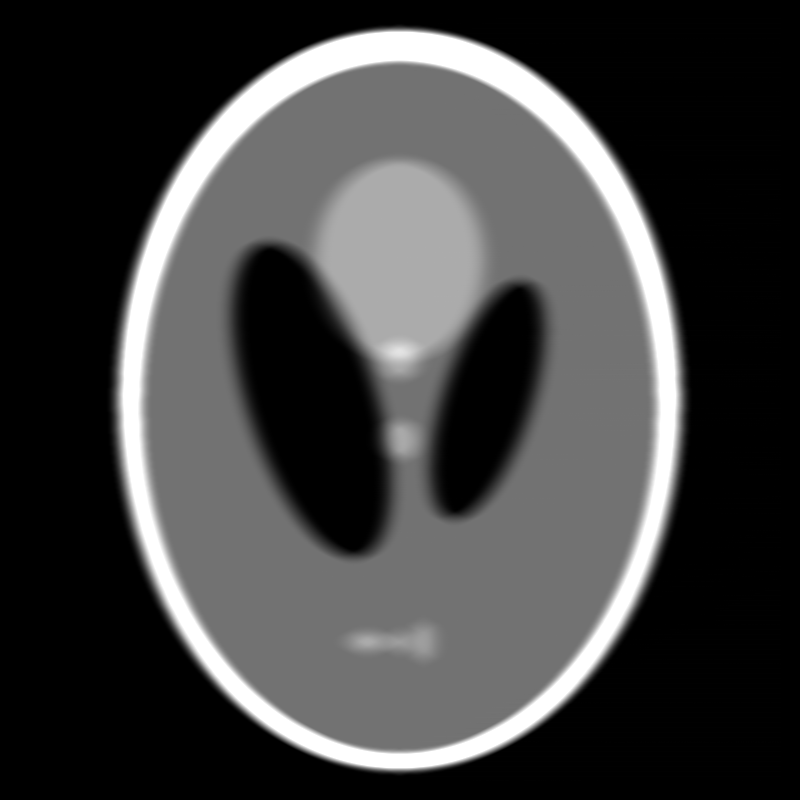}
		\caption{Effect of a big $\epsilon$ on reconstruction.}
		\label{fig:fig3}	
	\end{subfigure}
	\caption{Reconstruction using formula (\ref{VRT-inversion-2}).}
	\label{fig-rec}
	\end{center}		
\end{figure}

 The second reconstruction is based on formula (\ref{VRT-inversion-2}) and uses the averages of $f$ over infinitesimal parallelograms. Applying the averaging formula from Theorem \ref{ftc} and fixing $t=\epsilon$ as a constant representing the side length of the infinitesimal parallelogram, we get the following reconstructions (see Figure \ref{fig-rec}).

Note that in our numerical reconstruction we have chosen a specific size for the infinitesimal parallelogram. To get the best outcome we need to find an appropriate value for  $\epsilon$. Large $\epsilon$ will produce a more blurry outcome and a small value leads to more artifacts.

\subsection{Effects of noise}
When we have noise in the broken line data, we can refine the value of $\epsilon$ based on expected noise in the input. For  illustration  of the effects of Gaussian noise in our reconstructions see Figure \ref{fig:fig4}.

In practice, we  can also apply an averaging filter with appropriate window size on the broken line data to get a better reconstruction. See Figure \ref{fig:fig5}.

\begin{figure}[H]
	\begin{subfigure}{.5\textwidth}
		\centering
		\includegraphics[width=.7\linewidth]{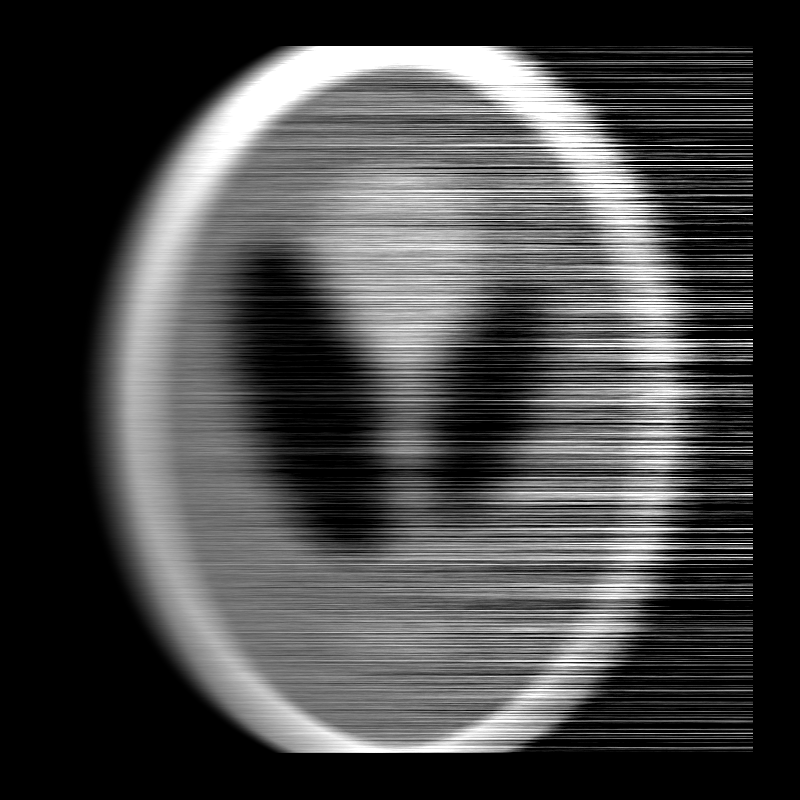}
		\caption{10 percent noise, $\epsilon=23$}
		\label{fig:sfig21}
	\end{subfigure}%
	\begin{subfigure}{.5\textwidth}
		\centering
		\includegraphics[width=.7\linewidth]{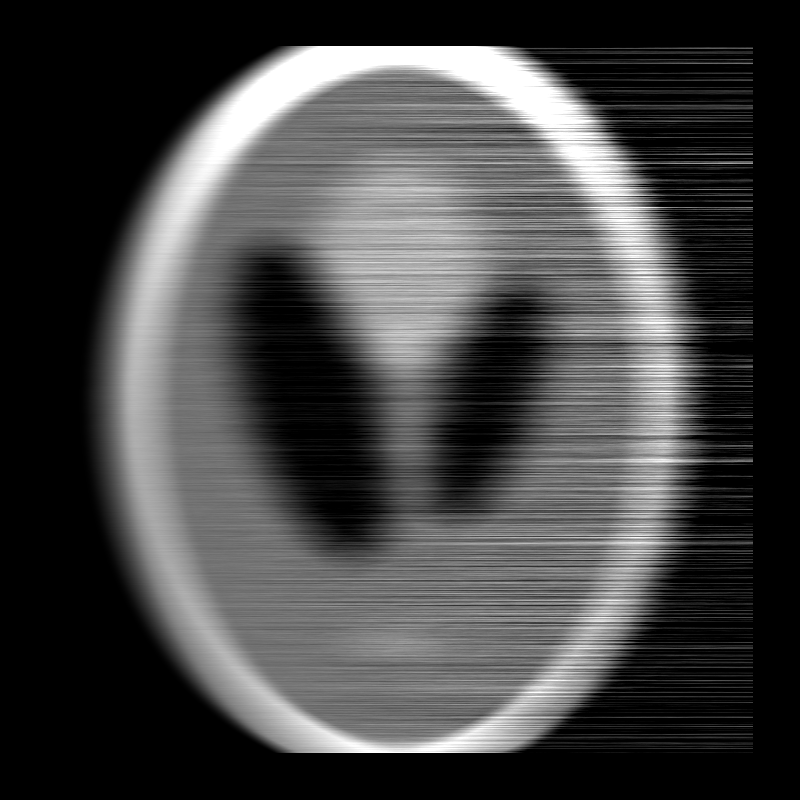}
		\caption{5 percent noise, $\epsilon=23$}
		\label{fig:sfig22}
	\end{subfigure}
\end{figure}
\begin{figure}[H]
\ContinuedFloat
	\begin{subfigure}{.5\textwidth}
		\centering
		\includegraphics[width=.7\linewidth]{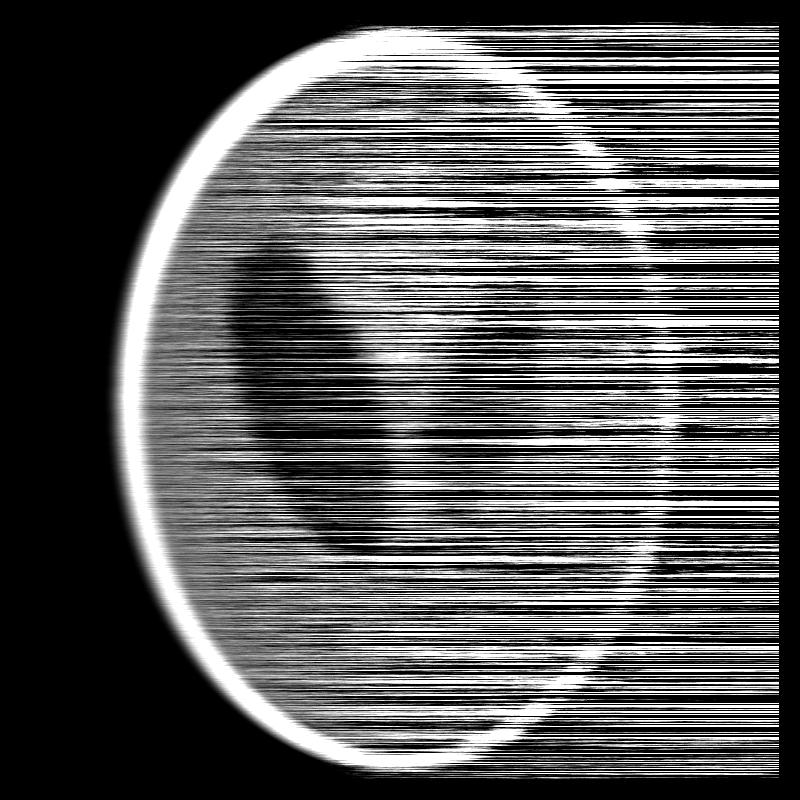}
		\caption{10 percent noise, $\epsilon=10$}
		\label{fig:sfig32}
	\end{subfigure}
	\begin{subfigure}{.5\textwidth}
		\centering
		\includegraphics[width=.7\linewidth]{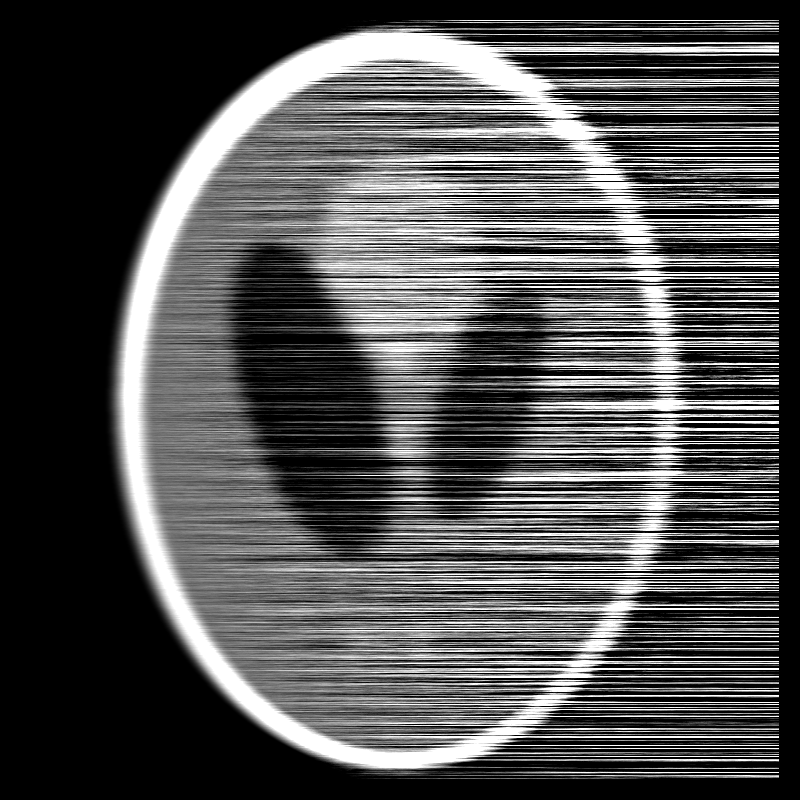}
		\caption{5 percent noise, $\epsilon=10$}
		\label{fig:sfig42}
	\end{subfigure}
	\caption{The effect of noise with different values of $\epsilon$ on reconstruction}
	\label{fig:fig4}
\end{figure}

\begin{figure}[H]
	\begin{center}
		\includegraphics[width=.35\linewidth]{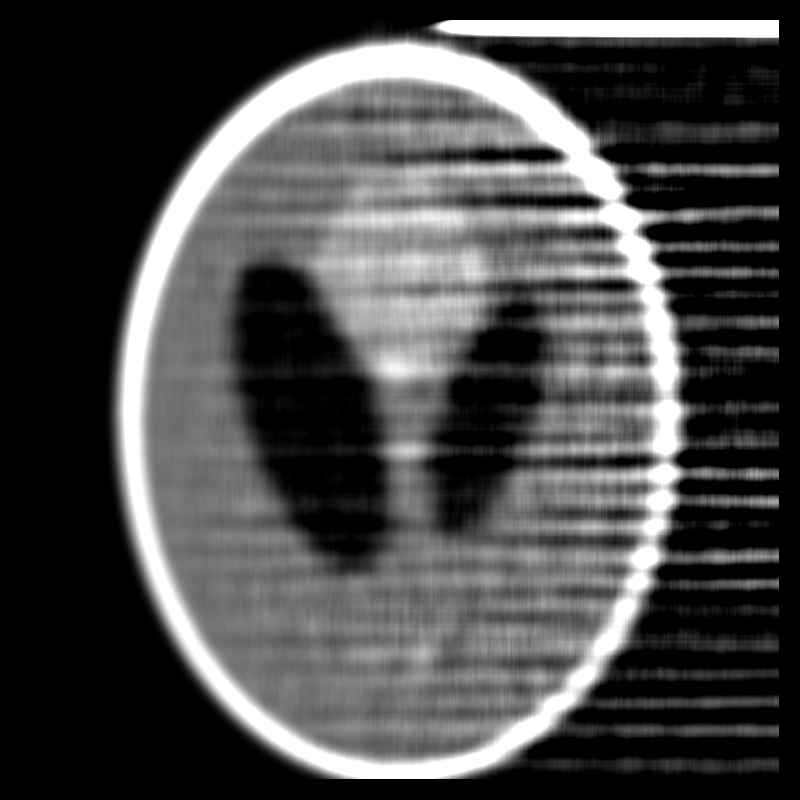}
		\caption{Filtering noisy VLT data. Window size = 12, $\epsilon=12$ and the noise is 10 percent.}
		\label{fig:fig5}
	\end{center}
	
\end{figure}


\section{Additional Remarks}\label{Chap_remarks}

\begin{enumerate}
\item {\bf Essence of Corners.} Our methodology has a very intuitive geometric interpretation, which in a nutshell can be described as follows. Use the Radon data to obtain a weighted average of the image function on a compact set (a polygon or a polyhedron), and then take the limit of that quantity (when the size of the set is sent to zero) to recover the function. The first step of that process can be accomplished relatively easily due to the presence of ``corners'' in the trajectories (surfaces) of integration, which distinguishes the VLT and the CRT from the conventional generalized Radon transforms integrating over smooth surfaces.
\item {\bf General Approach for Manifolds.} The methodology presented in this paper suggests a framework for deriving similar results on manifolds. Namely, one can use generalizations of the coarea formula in geometric measure theory to compute the corresponding $F$ as in Theorem \ref{inversion1}. Then properly combining values of $F$ at different points one can get (weighted) integrals of $f$ over ``nicely shrinking sets'', which may be used to produce $f$ through the Lebesgue differentiation theorem.

\item {\bf Weak Solutions.} Some of the prior work on inversion of the VLT and the CRT has been based on PDE techniques (e.g. see \cite{Kats_Krylov-13, Palamodov}). In a certain sense, our approach to solving these problems can be interpreted as finding weak solutions to the corresponding problems, i.e. satisfying the appropriate integral equations (our solutions are not necessarily differentiable).

\item {\bf A range description for the polyhedral case} may be derived with an appropriate (albeit very tedious) generalization of the techniques used in the case of the VLT.

\end{enumerate}

\section{Summary}\label{Chap_sum}

The paper presents a new approach to the inversion of a class of generalized Radon transforms, which map a function to its integrals along broken lines in the plane or polyhedral cones in higher dimensions.
These types of transformations play an important role in several modern imaging modalities based on physics of scattered particles.

We derived new explicit inversion formulae for the VLT and the CRT,  as well as  re-proved  some  previously  known  results  using more intuitive geometric ideas.  Using our inversion method for the VLT, we described the range of that transform when applied to a fairly broad class of functions, and proved some support theorems.  The efficiency of our  method  was  demonstrated  on  several  numerical  examples.  As an  auxiliary  result  that  played  a  big  role  in  this  article,  we  derived  a generalization  of  the  Fundamental  Theorem of  Calculus,  which  we called the Cone Differentiation Theorem.

\section{Appendix I: Absolute Continuity}\label{Chap_ap1}
In the case of finite measure $\nu$, Definition \ref{def-abscon} of absolute continuity is equivalent to the following (e.g. see Theorem 6.11 in \cite{Rudin:1987:RCA:26851}):

A measure $\nu$ is absolutely continuous with respect to Lebesgue measure $\mu$ if for any $\epsilon > 0$ there exists a $\delta > 0 $ such that $\mu(E) < \delta$ implies $\nu(E) < \epsilon$. \\

A cumulative distribution function $F$ is called absolutely continuous if for any $\epsilon > 0$ there exist a $\delta > 0 $ such that $\sum \mu(P_i) < \delta$ implies $\sum \nu_F(P_i) < \epsilon$ for any finite collection of disjoint parallelepipeds $\{ P_i \}$. Here, the induced measure $\nu_F$ on parallelepipeds is defined using the cumulative distribution function $F$. \\

We want to prove the following statement (e.g. see Proposition 26 in Section 20.3 of \cite{Royden}):

\vspace{2mm}
\noindent $F$ is absolutely continuous $\iff $  $\nu_F$ induced by $F$ is absolutely continuous with respect to the Lebesgue measure $\mu$
\\

 \noindent $(\Rightarrow )$  Let $\epsilon > 0$ be given. Then, by the assumption of absolute continuity of $F$ we can choose $\delta > 0$ such that $\sum \mu(P_i) < \delta$ implies $\sum \nu_F(P_i) < \epsilon/2$ for any finite collection of disjoint parallelepipeds $\{ P_i \}$. Now, let $E\subset \mathbb{R}^n$ with $\mu(E)< \delta/2$. Then there exists a countable disjoint collection of parallelepipeds $\{ P_i \}$ such that $E \subset \bigcup\limits_{i=1}^{\infty} P_i$ and $\mu\left(\bigcup\limits_{i=1}^\infty P_i\right)=\sum\limits_{i=1}^\infty \mu(P_i)< \delta$.  This follows from the definition of the outer measure, the Lebesgue cover and the fact that our parallelepipeds are half-open (see formula (\ref{parallele})).
 For the measure $\nu_F$ we have
 $$
 \nu_F\left(\bigcup\limits_{i=1}^{n} P_i\right)=\sum\limits_{i=1}^n \nu_F(P_i)<\epsilon/2,
 $$
 and hence
 $$
 \nu_F\left(\bigcup\limits_{i=1}^{\infty} P_i\right)=\lim\limits_{n\to\infty}\nu_F\left(\bigcup\limits_{i=1}^{n} P_i\right)=\lim\limits_{n\to\infty}\sum\limits_{i=1}^n \nu_F(P_i)\le\epsilon/2.
 $$
 Finally, we have
 $$
 \nu_F(E) \leq \nu_F\left(\bigcup\limits_{i=1}^{\infty} P_i\right) < \epsilon.
 $$

\vspace{3mm}

 \noindent $(\Leftarrow )$ Given $\epsilon > 0$ we pick $\delta>0$ using the definition of absolute continuity for $\nu_F$. Let $\{ P_i \}$ be any disjoint collection of parallelepipeds with $\sum \mu(P_i) < \delta$ and hence $\mu (\cup P_i) < \delta$. Then by the hypothesis $\sum \nu_F(P_i) = \nu_F (\cup P_i) < \epsilon$.

 \section{Appendix II: VLT as a Map on $L^1(\mathbb{R}^2)$}\label{Chap_ap2}

 A ray has measure zero as a subset of $\mathbb{R}^2$, hence changing the values of a function $f\in L^1(\mathbb{R}^2)$ along a ray will produce an equivalent function in $L^1$ sense.  Here we  show that the ray transform respects this equivalence relation, i.e. the ray transform maps two equivalent functions to the same equivalence class. The same argument will work for VLT, as it consists of a sum of two ray transforms along fixed directions. Without loss of generality we will prove the statement for the case of vertical rays.

 Assume $f,g\in L^1(\mathbb{R}^2)$ and $f=g$ almost everywhere. For an arbitrary $c\in\mathbb{R}$ consider the formal notations
 \[
 \phi_c(x) = \int_c^{\infty} f(x,y)\,dy,
 \]
 \[
 \psi_c(x) = \int_c^{\infty} g(x,y)\,dy.
 \]

By Fubini's theorem for any real numbers $a<b$ we have
\[
\int_{[a,b]\times [c,\infty]} f\, d\mu = \int_a^b \int_c^{\infty} f(x,y) \, dy\, dx= \int_a^b \phi_c(x)\,dx ,
\]
\[
\int_{[a,b]\times [c,\infty]} g\, d\mu = \int_a^b \int_c^{\infty} g(x,y) \, dy\, dx= \int_a^b \psi_c(x)\,dx,
\]
where $\phi_c(x)$ and $\psi_c(x)$ are integrable functions of $x$ and
\[
\int_a^b \phi_c(x)\,dx  = \int_a^b \psi_c(x)\,dx.
\]
By the Lebesgue Differentiation Theorem we get $\phi_c(x)=\psi_c(x)$ for almost every $x$.
Moreover, this statement is true for all $c\in\mathbb{R}$.

Now, let us recall some facts about product measure. Let $(X_1, \Sigma_1, \mu_1)$ and $(X_2, \Sigma_2, \mu_2)$ be  $\sigma$-finite measure spaces. Then the product measure on the product measurable space satisfies the following:
$$
(\mu_1 \times \mu_2)(E)= \int\limits_{X_2} \mu_1(E^y)\, d\mu_2(y),
$$
where $E^y=\left\{x\in X_1|\, (x,y)\in E\right\}$.

Hence, if for some $E\subseteq \mathbb{R}^2$ we have $(\mu_1 \times \mu_2)(E)\ne 0$, then $\mu_1(E^y)\ne 0$ for some $y$. Now if we let, $E=\{(x,c)\in \mathbb{R}^2 \, |\,  \phi_c(x) \neq \psi_c(x) \}$, this will imply that $E$ is of  measure zero.

In other words, we just showed, that if $f,g\in L^1(\mathbb{R}^2)$ and $f=g$ almost everywhere in $\mathbb{R}^2$, then the values of the ray transform of these two functions coincide for almost every (vertical) ray in the plane.

\section{Acknowledgements}

The work of both authors was partially supported by NSF grant DMS 1616564. G. Ambartsoumian was also supported in part by Simons Foundation grant 360357.

A part of this paper was written during G. Ambartsoumian's visit to the American University of Armenia (AUA). He thanks the administration and the staff of AUA for their support and stimulating environment.

\bibliographystyle{plain}
\bibliography{references}

\end{document}